\documentclass[leqno]{article}
\usepackage{amsmath,amssymb,amsthm,amscd,dsfont}
\numberwithin{equation}{section} \allowdisplaybreaks
\begin{document}
\newtheorem{theorem}{Theorem}[section]
\newtheorem{defin}{Definition}[section]
\newtheorem{prop}{Proposition}[section]
\newtheorem{corol}{Corollary}[section]
\newtheorem{lemma}{Lemma}[section]
\newtheorem{rem}{Remark}[section]
\newtheorem{example}{Example}[section]
\title{Foliated Lie and Courant Algebroids}
\author{{\small by}\vspace{2mm}\\Izu Vaisman}
\date{}
\maketitle
{\def\thefootnote{*}\footnotetext[1]%
{{\it 2000 Mathematics Subject Classification: 53C12, 53D17} .
}}
\begin{center} \begin{minipage}{12cm}
A{\footnotesize BSTRACT. If $A$ is a Lie algebroid over a foliated
manifold $(M,\mathcal{F})$, a foliation of $A$ is a Lie subalgebroid
$B$ with anchor image $T\mathcal{F}$ and such that $A/B$ is locally
equivalent with Lie algebroids over the slice manifolds of
$\mathcal{F}$. We give several examples and, for foliated Lie
algebroids, we discuss the following subjects: the dual Poisson
structure and Vaintrob's super-vector field, cohomology and
deformations of the foliation, integration to a Lie groupoid. In the
last section, we define a corresponding notion of a foliation of a
Courant algebroid $A$ as a bracket-closed, isotropic subbundle $B$
with anchor image $T\mathcal{F}$ and such that $B^\perp/B$ is
locally equivalent with Courant algebroids over the slice manifolds
of $\mathcal{F}$. Examples that motivate the definition are given.}
\end{minipage}
\end{center} \vspace*{5mm}
The main categories of interest for Differential Geometry are the
$C^\infty$ category and the complex analytic category. However,
there also is a significant interest in the category of foliated
manifolds. On the other hand, in the last thirty years Lie
algebroids became a central theme of differential-geometric
research, usually within the framework of the $C^\infty$ category.
Recently, a general study of holomorphic Lie algebroids has also
been done (see \cite{{GSX},{GSX2}}). The aim of the present paper is
to start a similar study of Lie algebroids in the
$C^\infty$-foliated category.

Essentially, we will say that a Lie subalgebroid $B$ is a
foliation of the Lie algebroid $A$ over the foliated manifold
$(M,\mathcal{F})$ if the anchor image of $B$ is $T\mathcal{F}$ and
$A/B$ is locally equivalent with Lie algebroids over the slice
manifolds of $\mathcal{F}$. In Section 1, after giving the precise
definitions and first properties, we discuss several examples of
foliated Lie algebroids $(A,B)$: classical foliations, transitive
foliations that correspond to foliated principal bundles, foliated
Dirac structures, etc. In Section 2, we show that the tangent Lie
algebroid of a foliated Lie algebroid is a foliated Lie algebroid.
In Section 3, we establish the characteristic properties of the
dual Poisson structure and of Vaintrob's super-vector field of a
foliated Lie algebroid. In Section 4, we define the cohomological
spectral sequence of a foliated Lie algebroid and prove a
Poincar\'e lemma in a particular case. Then, we define
deformations of a foliation $B$ of $A$ and prove the existence of
corresponding, cohomological, infinitesimal deformations. In
Section 5, we show that, if $A$ is integrable to a Lie groupoid
$G$ and $B$ is a foliation of $A$, $G$ has a foliation
$\mathcal{G}$ such that the restriction of the construction of the
Lie algebroid of $A$ to the leaves of $\mathcal{G}$ produces the
Lie subalgebroid $B$. Finally, in Section 6, we give a
corresponding definition of foliated Courant algebroids. With a
notation similar to the above, in the Courant case we will ask $B$
to be a bracket-closed, isotropic subbundle of the Courant
algebroid $A$, with anchor image $T\mathcal{F}$, such that
$B^\perp/B$ is locally equivalent with Courant algebroids over the
slice manifolds of $\mathcal{F}$. A usual foliation is a foliation
of the Courant algebroid $TM\oplus T^*M$ in this sense.

In the whole paper, we will use the Einstein summation convention.
\section{Definitions and examples}
We begin by recalling some basic facts concerning foliations; the
reader may find details in \cite{Mol}. Consider an $m$-dimensional
manifold $M$ endowed with a $C^\infty$-foliation $ \mathcal{F}$ with
$n$-dimensional leaves and codimension $k=m-n$. Then, each point has
local adapted coordinates $(x^a,y^u)$
$(a,b,...=1,...,k;\,u=1,...,n)$, with origin-centered, rectangular
range, such that the local equations of $ \mathcal{F}$ are
$x^a=const.$ and $y^u$ are leaf-wise coordinates. The foliation has
the following, associated exact sequence of vector bundles
\begin{equation}\label{sirfol} 0\rightarrow T\mathcal{F}
\stackrel{\iota}{\rightarrow}TM\stackrel{\psi}{\rightarrow}\nu\mathcal{F}
\rightarrow0,\end{equation} where $T\mathcal{F}$ is the tangent
bundle of the leaves, which will also be denoted by $F$, and
$\nu\mathcal{F}=TM/F$ is the transversal bundle of the foliation.

A mapping
$\varphi:(M_1,\mathcal{F}_1)\rightarrow(M_2,\mathcal{F}_2)$ between
two foliated manifolds is a {\it foliated mapping} if it sends
leaves of $ \mathcal{F}_1$ into leaves of $ \mathcal{F}_2$. The
local equations of such a map are of the form
\begin{equation}\label{aplfol} x_2^{a_2}=x_2^{a_2}(x_1^{a_1}),\;
y_2^{u_2}=y_2^{u_2}(x_1^{a_1},y_1^{u_1}).\end{equation} The
mapping $\varphi$ is called a {\it leaf-wise immersion,
submersion, local diffeomorphism} if the morphism induced by its
differential $\varphi_*$ between the tangent bundles
$T\mathcal{F}_1,T\mathcal{F}_2$ is an injection, surjection,
isomorphism, respectively. There also exists a similar notion of a
{\it transversal immersion, submersion, local diffeomorphism} that
refers to the transversal bundles
$\nu\mathcal{F}_1,\nu\mathcal{F}_2$.

Let $\pi_{P\rightarrow M}:P\rightarrow(M,\mathcal{F})$ (notice our
notation for projections of bundles onto the base manifold; the
index of $\pi$ will be omitted if no confusion is feared) be a
principal $G$-bundle over the foliated manifold $(M,\mathcal{F})$.
$P$ is a {\it foliated bundle} if it is endowed with a foliation $
\mathcal{F}^P$ that satisfies the conditions:

(i) $T\mathcal{F}^P$ has a trivial intersection with the tangent
spaces of the fibers (i.e., the intersection is equal to $0$ and
$T\mathcal{F}^P$ is ``horizontal"),

(ii) $\mathcal{F}^P$ is invariant by $G$-right translations,

(iii) the projection $\pi$ is a foliated leaf-wise submersion.

In fact, because of (ii) $\pi$ induces a covering map between the
leaves of $ \mathcal{F}^P$ and the corresponding leaves of $
\mathcal{F}$. Equivalently, a foliated principal bundle may be
characterized by an atlas of local trivializations with $G$-valued,
foliated transition functions where $G$ is foliated by points. A
foliated principal bundle $P$ may be seen as the glue-up of
pullbacks of $G$-principal bundles on the local transversal
manifolds of the leaves of $ \mathcal{F}$.

The associated bundles of a foliated principal bundle are called
foliated bundles too; they also have an induced foliation with
leaves that cover the leaves of $ \mathcal{F}$. Accordingly, a
vector bundle is foliated if the corresponding principal bundle of
frames is foliated. It follows that a foliated structure on a
vector bundle $\pi:A\rightarrow M$ is a maximal
local-trivialization atlas, defined on neighborhoods $\{U\}$,
where the local bases $(b_i)$ are related by matrices of
$\mathcal{F}$-foliated functions. Then, $A$ has the $
\mathcal{F}$-covering foliation $ \mathcal{F}^A$ and a cross
section $s$ of $A$ is foliated if the corresponding mapping
$s:(M,\mathcal{F})\rightarrow(A,\mathcal{F}^A)$ is foliated. Over
a trivializing neighborhood $U$ we have $s=\alpha^ib_i$ where
$\alpha^i$ are $\mathcal{F}$-foliated functions. Furthermore, we
have a sheaf $\Phi_{ \mathcal{F}}(A)$ of germs of foliated
cross-sections, which is a locally free sheaf-module over the
sheaf $\Phi_{\mathcal{F}}$ of germs of foliated functions on $M$
and which spans the sheaf $\Phi(A)$ of germs of cross sections of
$A$ over $C^\infty(M)$. Using the fact that $\Phi(A)$ is a locally
free sheaf it follows that every subsheaf $\Phi_{
\mathcal{F}}(A)\subseteq\Phi(A)$ that is a locally free sheaf-module over
$\Phi_{\mathcal{F}}$ and spans $\Phi(A)$ yields a well defined
foliated structure on $A$. Indeed, let $a_i$ be a basis of local
cross sections of $A$. Since $\Phi_{ \mathcal{F}}(A)$ spans
$\Phi(A)$, we have $a_i=\sum\beta_\alpha^\lambda b_\lambda$ with the
germs of $b_\lambda\in\Phi_\mathcal{F}(A)$. The total number of
elements $b_\lambda$ is finite and an independent subsystem may be
used as a foliated local basis. Two such bases are related by
foliated transition functions.

A foliated structure of $A$ may be identified with a family of
vector bundles $A_{Q_U}$ over quotient spaces
$Q_U=U/_{U\cap\mathcal{F}}$ such that $A|_U=\pi_{U\rightarrow
Q_U}^{-1}(A_{Q_U})$ ($\{U\}$ is an open covering of $M$ and
$\pi^{-1}$ denotes bundle pullback). Conversely, a family
$\{A_{Q_U}\}$ defines a foliated bundle $A$ if $\pi_{U\rightarrow
Q_U}^{-1}(A_{Q_U}|_{U\cap V})=\pi_{V\rightarrow
Q_V}^{-1}(A_{Q_V}|_{U\cap V})$. If $
\mathcal{F}$ consists of the fibers of a submersion $M\rightarrow
Q$, and if these fibers are assumed to be connected and simply
connected, then a foliated bundle $A$ over $M$ is the pullback of
a projected bundle $A_Q\rightarrow Q$ (see Lemma 2.5 and
Proposition 2.7 of \cite{Mol}).

Finally, a {\it foliated morphism} $\phi:E_1\rightarrow E_2$
between two foliated vector bundles over the same basis
$(M,\mathcal{F})$ is a vector bundle morphism that sends foliated
cross sections to foliated cross sections.

The most important example of a foliated vector bundle is the
transversal bundle $\nu\mathcal{F}$ of the foliation $\mathcal{F}$.
Notice that the foliated cross sections
$s\in\Gamma_{fol}\nu\mathcal{F}$ act on foliated functions $f\in
C^\infty_{fol}(M,\mathcal{F})$ (the index ``fol" is a shortcut for
``foliated") by $s(f)=Xf\in C^\infty_{fol}(M,\mathcal{F})$, for
every $X$ such that $s=[X]_{{\rm mod.}\, F}$ ($s(f)$ is well defined
since $X$ is a {\it foliated vector field}, i.e., a foliated cross
section $X:(M,\mathcal{F})\rightarrow(TM,T\mathcal{F})$).
Furthermore, on the space $\Gamma_{fol}\nu\mathcal{F}$, there exists
a well defined Lie algebra bracket
\begin{equation}\label{crosetpenu}[[X_1]_{{\rm mod.}\,F},
[X_2]_{{\rm mod.}\,F}]_{\nu\mathcal{F}}=[X_1,X_2]_{{\rm
mod.}\,F}\end{equation} induced by the Lie bracket of the
corresponding foliated vector fields.

Accordingly, the definition of a Lie algebroid suggests the
following definition.
\begin{defin}\label{defalgtrans} {\rm An
$\mathcal{F}$-{\it transversal-Lie algebroid} over
$(M,\mathcal{F})$ is a foliated bundle $E$ endowed with a Lie
algebra bracket $[\,,\,]_E$ on $\Gamma_{fol}E$ and a foliated,
{\it transversal anchor} morphism $\sharp_E:E\rightarrow\nu
\mathcal{F}$ such that $$\begin{array}{l}
1)\;\sharp_E[e_1,e_2]_E=[\sharp_E e_1,\sharp_E
e_2]_{\nu\mathcal{F}},\;\forall
e_1,e_2\in\Gamma_{fol}E,\vspace{2mm}\\
2)\;[e_1,fe_2]_E=f[e_1,e_2]_E+(\sharp_E e_1(f))e_2,\;\forall
e_1,e_2\in\Gamma_{fol}E,f\in
C^\infty_{fol}(M).\end{array}$$}\end{defin}

The symbols $\sharp,\,[\,,\,]$ will be used for any
transversal-Lie and Lie algebroid while the index that denotes the
bundle will be omitted if there is no risk of confusion. Notice
also that we may not ask the morphism $\sharp_E$ to be foliated, a
priori; this follows from condition 2), which implies $\sharp
e_1(f)\in C^\infty_{fol}(M,\mathcal{F})$, $\forall f\in
C^\infty_{fol}(M,\mathcal{F})$. The following result is obvious
\begin{prop}\label{obsproiectii} If $E$ is a transversal-Lie algebroid
over $(M,\mathcal{F})$, $\nu\mathcal{F}$ projects to the tangent
bundles of the local quotient spaces $Q_U$ of $\mathcal{F}$ and
the local projected bundles $\{E_{Q_U}\}$ are Lie algebroids. In
the case of a submersion $M\rightarrow Q$ with connected and
simply connected fibers, $E\rightarrow M$ projects to a Lie
algebroid over $Q$.\end{prop}

A transversal-Lie algebroid is not a Lie algebroid, and we shall
define the notion of a foliated Lie algebroid as follows.

If $A$ is a Lie algebroid over $(M,\mathcal{F})$ and if $B$ is a
Lie subalgebroid, a cross section $a\in\Gamma A$ will be called a
$B$-{\it foliated cross section} if $[b,a]_A\in\Gamma B$, $\forall
b\in\Gamma B$. The correctness of this definition follows from the
fact that $[b,a]_A\in\Gamma B$ implies
$[fb,a]_A=f[b,a]_A-(\sharp_A a(f))b\in\Gamma B$. We shall denote
by $\Gamma_BA$ the space of $B$-foliated cross sections of $A$.
The Jacobi identity shows that the space $\Gamma_BA$ is closed by
$A$-brackets. If $\sharp_A(B)\subseteq T\mathcal{F}$ then,
$\forall a\in\Gamma_BA$ and $\forall f\in
C^\infty_{fol}(M,\mathcal{F})$, we have $[b,fa]_A\in\Gamma_BA$.
Now, we give the following definition
\begin{defin}\label{folLiealg}
{\rm If $A$ is a Lie algebroid on $(M,\mathcal{F})$, a {\it
foliation of $A$ over $ \mathcal{F}$} is a Lie subalgebroid
$B\subseteq A$ that satisfies the following two conditions: 1)
$\sharp_A(B)=F$, 2) locally, $\Gamma A$ is spanned by $\Gamma_BA$
over $C^\infty(M)$ (this is called the {\it foliated-generation
condition}). A pair $(A,B)$ as above will be called a {\it foliated
Lie algebroid}. If $\sharp_A:B\rightarrow F$ is a vector bundle
isomorphism, $(A,B)$ is a {\it minimally foliated Lie
algebroid}.}\end{defin}
\begin{rem}\label{Frobenius} {\rm The condition that $B$ is
closed by $A$-brackets may be replaced by the Frobenius condition:
$\forall\alpha\in ann\,B$, the $A$-exterior differential
$d_A\alpha$ belongs to the ideal of $A$-forms generated by
$ann\,B$.}\end{rem}
\begin{lemma}\label{bazefol} If $B$ is a foliation of $A$ and $C$
is a complementary subbundle of $B$ in $A$ ($A=B\oplus C$),
$\Gamma C$ has local bases that consist of $B$-foliated cross
sections of $A$.\end{lemma}
\begin{proof} Let $(c_h)$ be a local basis of $\Gamma C$ over an
open neighborhood $U\subseteq M$. By the foliated-generation
condition we have
$$c_h=\gamma_h^ua_{(h)u}=\gamma_h^upr_Ca_{(h)u},$$ where
$\gamma_h^u$ are differentiable functions and $a_{(h)u}$,
therefore, $pr_Ca_{(h)u}$ too, are local $B$-foliated cross
sections of $A$. Since the set $\{pr_Ca_{(h)u}\}$ is finite, after
shrinking the neighborhood $U$ as necessary, there exists a subset
of linearly independent, local cross sections $\{pr_Ca_v\}$ such
that $c_h=\lambda_h^va_{v}pr_Ca_v$. This subset obviously is the
required local basis of $\Gamma C$.
\end{proof}
\begin{rem}\label{obsCF} {\rm
The complementary subbundle $C$ may have elements with projection
in $F$ since we did not ask $B=\sharp_A^{-1}(F)$. On the other
hand, by condition 1) of Definition \ref{folLiealg}, there exist
decompositions $B=ker\,\sharp_B\oplus P$ where  $P$ is a subbundle
of $B$.}\end{rem}
\begin{prop}\label{ApeB} If $B$ is a foliation of the Lie
algebroid $A$, then the vector bundle $E=A/B$ is a transversal-Lie
algebroid on $ (M,\mathcal{F})$.\end{prop}
\begin{proof} Consider a decomposition $A=B\oplus C$ and take
$B$-foliated local bases of $C$ over an open covering $\{U\}$ of
$M$, which exist by Lemma \ref{bazefol}. Let
$\tilde{a}_v=\lambda_v^w a_w$ be a transition between two such bases
and take the bracket by $b\in\Gamma B$. We get $(\sharp
b)\lambda_u^v=0$ and, since $\sharp(B)=F$, the functions
$\lambda_u^v$ must be $ \mathcal{F}$-foliated. Thus, $C$ has an
induced foliated structure, which transfers to $A/B\approx C$. We
notice that $[a]_{{\rm mod.}\,B}\in\Gamma_{fol}(A/B)$ if and only if
$a\in\Gamma_BA$. Furthermore, the bracket of $B$-foliated cross
sections of $A$, which is $B$-foliated too, descends to a well
defined Lie bracket on $\Gamma_{fol}E$. Finally, the anchor
$\sharp_A$ induces an anchor $\sharp_E:E\rightarrow\nu\mathcal{F}$
as required by Definition \ref{defalgtrans}. The fact that
$\sharp_E$ is a foliated morphism may be justified as in the
observation that follows Definition \ref{defalgtrans}.
Alternatively, we can see that, $\forall a\in\Gamma_BA$, $\sharp_Aa$
is a foliated vector field on $(M,\mathcal{F})$  as follows:
condition 1) of Definition \ref{folLiealg} shows that a vector field
$Y\in\Gamma F$ is of the form $Y=\sharp_Ab$, $b\in\Gamma B$, hence,
$$[Y,\sharp_Aa]_A=[\sharp_Ab,\sharp_Aa]_A=\sharp_A[b,a]\in\Gamma
F,$$ which is the required property.
\end{proof}
\begin{corol}\label{corolextensie}
A foliated Lie algebroid $(A,B)$ over $(M,\mathcal{F})$ produces a
commutative diagram
\begin{equation}\label{extension} \begin{array}{ccrcrcrcc}
0&\rightarrow& B&
\stackrel{j}{\rightarrow}&A&\stackrel{\sigma}{\rightarrow}&E&
\rightarrow&0\vspace{2mm}\\ & & \sharp_B\downarrow& &
\sharp_A\downarrow& & \sharp_E\downarrow& &\vspace{2mm}\\
0&\rightarrow& T\mathcal{F}&
\stackrel{\iota}{\rightarrow}&TM&\stackrel{\psi}{\rightarrow}&\nu\mathcal{F}
&\rightarrow&0\end{array}\end{equation} where the lines are exact
sequences of vector bundles.\end{corol}
\begin{defin}\label{defextensie} {\rm
If the transversal-Lie algebroid $E$ is foliately equivalent with
a quotient $A/B$ where $(A,B)$ is a foliated Lie algebroid, then
$E$ may be placed in a diagram (\ref{extension}) and $A$ will be
called an {\it extension} of $E$. An extension of $E$ such that
$B\approx F$ will be called a {\it minimal extension}.}\end{defin}

In what follows we give several examples of foliated Lie
algebroids.
\begin{example}\label{folclasica} {\rm For any foliation
$\mathcal{F}$ of a manifold $M$, $T\mathcal{F}$ is a foliation of
the Lie algebroid $TM$.}\end{example}
\begin{example}\label{subreg} {\rm If a Lie algebroid $A$ has a
regular subalgebroid $B$, $M$ has the foliation $\mathcal{F}$ such
that $F=\sharp(B)$. Then, $(A,B)$ is a foliated Lie algebroid if
the foliated-generation condition is satisfied. For instance, if
$A$ is regular, we may take $B=A$ and the foliated-generation
condition is trivially satisfied. Another obvious example is that
of a regular subalgebroid $B$ such that $B\subseteq
ker\,\sharp_A$. Then $B$ is a foliation of $A$ over the foliation
of $M$ by points. For instance, let $J^1A$ be the first jet Lie
algebroid of $A$ defined in \cite{CF}. It is well known that one
has the exact sequence
\begin{equation}\label{sirjeturi}
0\rightarrow Hom(TM,A)\stackrel{\subseteq}{\rightarrow}
J^1A\stackrel{\pi}{\rightarrow} A\rightarrow0,\end{equation} where
$\pi(j^1_xa)=a(x)$ $(x\in M,a\in\Gamma A)$ is a Lie algebroid
morphism, therefore, $ker\,\pi=Hom(TM,A)$ is a regular Lie
subalgebroid of $J^1A$. This subalgebroid is included in
$ker\,\sharp_{J^1A}$ because $j^1_xa\in ker\,\pi_x$ implies
$\sharp_{J^1A}(j_x^1a)=\sharp_A(a(x))=0$ (the equality
$\sharp_{J^1A}(j_x^1a)=\sharp_Aa(x)$ is a part of the definition
given in \cite{CF}). Thus, $Hom(TM,A)$ is a foliation of $J^1A$
over the points of $M$.}\end{example}
\begin{example}\label{Atiyah} {\rm The typical example of a
transitive Lie algebroid is $A=TP/G$, where $\pi:P\rightarrow M$
is a principal bundle of structure group $G$ over $M$ and the
quotient is the space of the orbits of $TP$ by the action of the
differential of the right action of $G$ on $P$. The cross sections
of $A$ may be identified with right-invariant vector fields on
$P$, which allows to define the bracket, and the anchor is induced
by $\pi_*:TP\rightarrow TM$ \cite{M}. Now, assume that $M$ has the
foliation $\mathcal{F}$ and $P$ is a foliated principal bundle
with the covering foliation $\mathcal{F}^P$ of $\mathcal{F}$.
Then, we get a Lie subalgebroid $B=T\mathcal{F}^P/G\subseteq A$
such that the restriction of the anchor of $A$ to $B$ has the
image $T\mathcal{F}$. The foliated-generation condition of $A$
with respect to $B$ is satisfied because it is satisfied for $TP$
with respect to $T\mathcal{F}^P$ and
$A/B=(TP/T\hat{\mathcal{F}})/G$. Therefore, $(A,B)$ is a foliated
Lie algebroid.}\end{example}
\begin{example}\label{vectfol} {\rm Following \cite{M},
if $P$ is the bundle of frames of a vector bundle $\pi:E\rightarrow
M$, Example \ref{Atiyah} has the following equivalent form. Let $A$
be the transitive Lie algebroid $\mathcal{D}(E)$ whose cross
sections are the linear differential operators $\Gamma
E\rightarrow\Gamma E$ of order $\leq1$. The bracket is the commutant
of the operators and the anchor is defined by the symbol of the
operator. In the presence of a foliation $\mathcal{F}$ on $M$ and if
we assume $E$ to be foliated, $\mathcal{D}(E)$ has a Lie
subalgebroid $B=\mathcal{D}_{\mathcal{F}}(E)$ defined by the
operators whose symbol is a vector field tangent to $\mathcal{F}$
($\mathcal{D}_{\mathcal{F}}(E)$ is commutant closed since
$F=T\mathcal{F}$ is closed by Lie brackets), and the anchor sends
$B$ onto the tangent bundle $F$: if $Y\in\Gamma F$, we get an
operator $D$ of symbol $Y$ by putting $Ds=0$ for
$s\in\Gamma_{fol}(E)$ and
\begin{equation}\label{Dt}Dt=\sum (Yf_i)s_i,\;\;\forall t=\sum
f_is_i,\;s_i\in\Gamma_{fol}E,f_i\in C^\infty(M)\end{equation} (the
definition is correct since if $s\in\Gamma_{fol}(E)$ and $f\in
C^\infty_{fol}(M,\mathcal{F})$ (\ref{Dt}) gives $D(fs)=0$). The
foliated-generation condition is satisfied too. Indeed, an operator
of symbol $X$ is of the form $D=\nabla_X+\phi$, where $\nabla$ is a
covariant derivative and $\phi\in End(E)$, therefore, $\phi\in
\mathcal{D}_{ \mathcal{F}}(E)$. The operator $D$ is foliated with
respect to $ \mathcal{D}_{ \mathcal{F}}(E)$ if and only if $X$ is
foliated with respect to $\mathcal{F}$ and it is known that the
vector fields on $M$ satisfy the foliated-generation condition.
Thus, the pair $(A,B)$ considered above is a foliated Lie
algebroid.}\end{example}
\begin{example}\label{dirac1} {\rm Let $D$ be a  Dirac structure on
$(M,\mathcal{F})$ such that $F\subseteq D$ ($D$ is
$\mathcal{F}$-projectable \cite{V1}). Then, $(D,F)$ is a foliated
Lie algebroid. The fact that the foliated-generation condition is
satisfied was proven in \cite{V1} (where, also, concrete examples
were given). We notice that there are foliated Dirac structures
where $F$ is strictly included in $\sharp_D^{-1}(F)$
$(\sharp_D=pr_{TM})$. For instance,
$$D=F\oplus\{(\sharp_P\alpha,\alpha)\,/\,
\alpha\in ann\,F\},$$ where $P$ is a foliated bivector field on $M$
such that the Schouten-Nijenhuis bracket $[P,P]$ vanishes on
$ann\,F$ (\cite{V1}, Example 5.1).}\end{example}
\begin{example}\label{dirac2} {\rm Let $D$ be a regular Dirac
structure on $M$ with the (regular) characteristic foliation
$\mathcal{E}$, $E=T\mathcal{E}=pr_{TM}D$ and let $\mathcal{F}$ be a
subfoliation of $\mathcal{E}$. Put
$$D_{ \mathcal{F}}=D\cap(F\oplus T^*M).$$ Notice that
the mapping $\psi:D/D_{\mathcal{F}}(x)\rightarrow E/F(x)$ given by
$[(X,\xi)]_{{\rm mod.}\,D_{F}}\mapsto [X]_{{\rm mod.}\,F}$ is an
isomorphism of vector spaces $\forall x\in M$. The existence of
this isomorphism shows that $dim\,D_{\mathcal{F}}(x)=const.$
Hence, $D_{\mathcal{F}}$ is a vector subbundle of $D$, which
obviously is a Lie subalgebroid of $D$ with anchor image $F$.
Moreover, since $\mathcal{F}$ is a subfoliation of $ \mathcal{E}$,
the vector bundle $E/F$ is $
\mathcal{F}$-foliated and so is its isomorphic image $D/D_{
\mathcal{F}}$. This shows that the pair $(D,D_{
\mathcal{F}})$ satisfies the foliated-generation condition and is
a foliated Lie algebroid.}\end{example}
\begin{prop}\label{extmin} For any transversal-Lie algebroid $E$
over the foliated manifold $(M,\mathcal{F})$ and any lift
$\rho:E\rightarrow TM$ of $\sharp_E$
$(\psi\circ\rho=\sharp_E,\,\psi=pr_{\nu\mathcal{F}})$ there exists
a canonical minimal extension $(A_0,\sharp_0,[\,,\,]_0)$ with
$A_0=F\oplus E$, $\sharp_0=Id+\rho$ and
\begin{equation}\label{crcanonic} \begin{array}{l}
[Y_1,Y_2]_0=[Y_1,Y_2],\;[Y,e]_0=[Y,\rho e],\;[e,Y]_0=[\rho
e,Y]\vspace{2mm}\\

[e_1,e_2]_0=([\rho e_1,\rho e_2]-\rho[e_1,e_2]_E)+[e_1,e_2]_E,
\end{array}\end{equation} $\forall Y,Y_1,Y_2\in\Gamma F$, $\forall
e,e_1,e_2\in\Gamma_{fol}E$, where the non-indexed brackets are Lie
brackets of vector fields. Conversely, for any minimal extension
$A$ of $E$, there exists a lift $\rho$ such that $A$ is isomorphic
to a twisted form of the canonical extension of $E$ by $F$.
\end{prop}
\begin{proof} The qualification ``twisted" is similar to that used
in the notion of a twisted Dirac structure and its exact meaning
will appear at the end of the proof. We shall use the content and
notation of diagram (\ref{extension}), which includes the mapping
$\psi=pr_{\nu\mathcal{F}}$. Notice that a lift $\rho$ has the
property that $\rho e$ is a foliated vector field $\forall
e\in\Gamma_{fol}E$. If we extend (\ref{crcanonic}) by
$$[Y,fe]_0=f[Y,\rho e]+(Yf)e,\;[e_1,fe_2]_0=f[e_1,e_2]_0+\rho
e_1(f)e_2,$$ $\forall f\in C^\infty(M)$, we get a skew-symmetric
bracket on $\Gamma A_0$. It is easy to see that the extension is
compatible with (\ref{crcanonic}) for $f\in
C^\infty_{fol}(M,\mathcal{F})$. Also, straightforward calculations
show that the axioms of a Lie algebroid, including the Jacobi
identity, are satisfied on arguments $Y\in\Gamma F,
e\in\Gamma_{fol}E$. Since these types of cross sections locally
generate $\Gamma A_0$, this proves the existence of the canonical
extension. Notice that, for $E=\nu\mathcal{F}$, $\rho$ is a
splitting of the lower line of (\ref{extension}) and, if we identify
$F\oplus\approx F\oplus im\,\rho=TM$, the previous construction
yields the standard Lie algebroid structure of $TM$.

For the converse, we start with diagram (\ref{extension}), which
yields $A=\phi(F)\oplus\tau(E)$ where $\phi:F\approx
B,\phi^{-1}=\sharp_B$ ($\phi$ exists because the extension is
minimal) and $\tau:E\rightarrow A$ is a splitting of the upper exact
sequence of the diagram ($\sigma\circ\tau=Id$). Then, we consider
the lift $\rho=\sharp_A\circ\tau$, which implies
$\sharp_A=\phi^{-1}+\rho\circ\sigma$. As for the brackets, we must
have
$$\phi^{-1}[\phi Y_1,\phi Y_2]_A=\sharp_A[\phi Y_1,\phi Y_2]_A
=[\sharp_A\phi Y_1,\sharp_A\phi Y_2]=[Y_1,Y_2],$$ i.e., $[\phi
Y_1,\phi Y_2]_A=\phi [Y_1,Y_2]$. Then, since $\sigma[\phi Y,\tau
e]_A=[\sigma\phi Y,\sigma\tau e]_E=0$, $[\phi Y,\tau e]_A\in\Gamma
B$ and
$$\phi^{-1}[\phi Y,\tau e]_A=\sharp_A[\phi Y,\tau e]_A
=[Y,\sharp_A\tau e]\,\Leftrightarrow\, [\phi Y,\tau
e]_A=\phi[Y,\rho e],\;\forall e\in\Gamma_{fol}E.$$ Finally, for
$e_1,e_2\in\Gamma_{fol}E$, we have
$$\sharp_A[\tau e_1,\tau e_2]_A-\rho[e_1,e_2]_E= [\rho e_1,\rho
e_2]-\rho[e_1,e_2]_E\in\Gamma F$$ because, if we apply $\psi$ to
the right hand side of the equality, we get
$$\psi[\rho e_1,\rho e_2]-\sharp_E[e_1,e_2]_E=[\psi\rho e_1,\psi\rho
e_2]_{\nu\mathcal{F}}-\sharp_E[e_1,e_2]_E$$ $$=[\sharp_E
e_1,\sharp_E e_2]_{\nu\mathcal{F}}-\sharp_E[e_1,e_2]_E=0$$ by the
axioms of a transversal-Lie algebroid. Now, we notice that, if $a\in
A$, $\sharp_Aa\in\Gamma F$ if and only if $a=\phi Y+\tau e$ where
$e\in ker\,\sharp_E$; indeed, by the commutativity of diagram
(\ref{extension}), $\psi\sharp_A\tau(e)=0$ is equivalent to
$\sharp_Ee=0$. This fact and the previous observation justify the
formula
\begin{equation}\label{twisted}[\tau e_1,\tau e_2]_A=([\tau e_1,\tau
e_2]_A-\tau[e_1,e_2]_E)+
\tau[e_1,e_2]_E\end{equation} $$=[\phi([\rho e_1,\rho
e_2]-\rho[e_1,e_2]_E)+\tau(\lambda(e_1,e_2))]+\tau[e_1,e_2]_E,$$
where $\lambda:\wedge^2E\rightarrow E$ is a $2$-form with values
in the kernel of $\sharp_E$. Therefore, $A$ is isomorphic by
$\phi+\tau$ to the canonical extension of $E$ by $F$ twisted by
the addition of the form $\lambda$. The latter must also be asked
to satisfy conditions that ensure the Jacobi identity.
\end{proof}
\begin{rem}\label{obsCarinena} {\rm In \cite{Car} the authors
define a reduction of Lie algebroids that, essentially, is a
projection along the fibers of a surjective submersion
$\pi:M\rightarrow M'$. Namely, $A\rightarrow M$ reduces to
$A'\rightarrow M'$ if there exists a Lie algebroid morphism
\cite{M} $(\Pi:A\rightarrow A',\pi:M\rightarrow M')$. One can see
that, if $B$ is a foliation of $A$ over the foliation of $M$ by
the fibers of $\pi$ and if the fibers are connected and simply
connected, then $A$ reduces to $A'$ defined by $A/B=\pi^{-1}A'$.}
\end{rem}
\section{The tangent Lie algebroid}
In this section we will show that the tangent Lie algebroid of a
foliated Lie algebroid $(A,B)$ over a foliated manifold
$(M,\mathcal{F})$ has a natural structure of a foliated Lie
algebroid over $(TM,T\mathcal{F})$. For this purpose, we formulate
the definition of the tangent Lie algebroid of a Lie algebroid
\cite{MX} in a suitable (but, not new) form.

Let $\pi:A\rightarrow M$ be an arbitrary vector bundle. Take a
neighborhood $U\subseteq M$ with local coordinates $(x^i)_{i=1}^m$,
with the local basis $(a_\alpha)^r_{\alpha=1}$ of $\Gamma A$
($r=rank A$) and the corresponding fiber coordinates $(\xi^\alpha)$.
On the intersection of two such neighborhoods $U\cap\tilde{U}$,
these coordinates have transition functions of the following local
form
\begin{equation}\label{transA}
\tilde{x}^i=\tilde{x}^i(x^j),\;\tilde{\xi}^\alpha=\Theta^\alpha_\beta(x^j)\xi^\beta\;\;
(\tilde{a}_\alpha=\Psi_\alpha^\gamma
a_\gamma,\;\Psi_\alpha^\gamma\Theta_\gamma^\beta=\delta_\alpha^\beta).\end{equation}

Correspondingly, one has natural coordinates $(x^i,\dot{x}^i)$ on
$T_UM$ and $(x^i,\xi^\alpha,$ $\dot{x}^i,\dot{\xi}^\alpha)$ on
$T_{\pi^{-1}U}A$, with the following transition functions:
\begin{equation}\label{transTA} \tilde{x}^i=\tilde{x}^i(x^j),\;
\tilde{\xi}^\alpha=\Theta^\alpha_\beta(x^j)\xi^\beta,\; \dot{\tilde{x}}^i=
\frac{\partial{\tilde{x}^i}}{\partial{{x}^j}}\dot{x}^j,\;
\dot{\tilde{\xi}}^\alpha=\frac{\partial\Theta^\alpha_\beta}{\partial
x^j}\xi^\beta\dot{x}^j
+\Theta^\alpha_\beta\dot{\xi}^\beta.\end{equation}

Formulas (\ref{transTA}) show the existence of the double vector
bundle
\begin{equation}\label{doublevb} \begin{array}{ccc}
TA&\stackrel{\pi_*}{\rightarrow}&TM\vspace{2mm}\\ \downarrow& &
\downarrow\vspace{2mm}\\
A&\stackrel{\pi}{\rightarrow}&M\end{array}\end{equation} where $TA$
is a vector bundle over $A$ with base coordinates $(x^i,\xi^\alpha)$
and fiber coordinates $(\dot{x}^i,\dot{\xi}^\alpha)$ and $TA$ is a
vector bundle over $TM$ with base coordinates $(x^i,\dot{x}^i)$ and
fiber coordinates $(\xi^\alpha,\dot{\xi}^\alpha)$. (Notice that the
cross sections of $\pi_*$ are not vector fields on $A$!)

If $A=TM$, the indices $\alpha,\beta,...$ may be replaced by
$i,j,...$ and $\Theta^i_j=
\partial\tilde{x}^i/\partial{x}^j$, which shows the
existence of the flip diffeomorphism $\phi:TTM\rightarrow TTM$
defined by the local coordinate equations
\begin{equation}\label{flip} \phi(x^i,\xi^i,\dot{x}^i,\dot{\xi}^i)
=(x^i,\dot{x}^i,\xi^i,\dot{\xi}^i).\end{equation}

On the fibration $TA\rightarrow A$ the fiber coordinates
$(\dot{x}^i,\dot{\xi}^\alpha)$ are produced by the local bases of
cross sections $(\partial/\partial
x^i,\partial/\partial\xi^\alpha)$. On the fibration
$\pi_*:TA\rightarrow TM$ the fiber coordinates
$(\xi^\alpha,\dot{\xi}^\alpha)$ are produced by local bases
$(c_\alpha,\partial/\partial\xi^\alpha)$, where the vector
$c_\alpha$ is the image of $a_\alpha$ under the natural
identification of the tangent space of the fibers of $A$ with the
fibers themselves, which are vector spaces. (The previous
assertion is justified by checking the invariance of the
expression $\xi^\alpha
c_\alpha+\dot{\xi}^\alpha(\partial/\partial\xi^\alpha)$ by the
coordinate transformations (\ref{transTA}).) It follows that for
$z\in TA$ the intersection space of the fibers of the two vector
bundle structures of $TA$ is
\begin{equation}\label{genP}
\mathcal{P}_z=\pi_{TA\rightarrow A}^{-1}(\pi_{TA\rightarrow
A}(z))\cap(\pi_{A\rightarrow M})^{-1}_*((\pi_{A\rightarrow
M})_*(z))=span\left\{\frac{\partial}{\partial\xi^\alpha}\right\}.\end{equation}
Therefore, $P=\cup_{z\in TA}\mathcal{P}_z$ is a vector subbundle
of $TA$, which we call the {\it bi-vertical subbundle}.

Furthermore, like for the tangent bundle $A=TM$, there are two
lifting processes of cross sections of $\pi:A\rightarrow M$ to
cross sections of $\pi_*:TA\rightarrow TM$.

One is the {\it complete lift} $a\mapsto a^C=a_*$, which sends the
cross section $a$ to the differential of the mapping $a:M\rightarrow
A$. The expression of $a_*$ by means of the local coordinates
$(x^i,\dot{x}^i)$ on $TM$ and
$(x^i,\dot{x}^i,\xi^\alpha,\dot{\xi}^\alpha)$ on $TA$ is
$$(x^i,\dot{x}^i)\mapsto(x^i,\dot{x}^i,\xi^\alpha(x^i),\dot{x}^i
\frac{\partial\xi^\alpha}{\partial x^i}),$$ where
$a=\xi^\alpha(x^i)a_\alpha$, whence we get
\begin{equation}\label{liftC} a^C=
\xi^\alpha c_\alpha + \dot{x}^i
\frac{\partial\xi^\alpha}{\partial x^i}
\frac{\partial}{\partial\xi^\alpha}.\end{equation}

The second lift, called the {\it vertical lift} and denoted by an
upper index $V$, is the isomorphism between the pullback
$\pi^{-1}_{TA\rightarrow A\rightarrow M}A$  of the vector bundle
$A$ to $TM$ and the bi-vertical subbundle $P$ defined by
\begin{equation}\label{liftV} a=\eta^\alpha a_\alpha\mapsto
a^V=\eta^\alpha\frac{\partial}{\partial\xi^\alpha}.\end{equation}
(Formulas (\ref{transA}) show that the bases $(a_\alpha)$ and
$(\partial/\partial\xi^\alpha)$ have the same transition
functions.)

Notice the following interpretation of the local basis
$(c_\alpha,\partial/\partial\xi^\alpha)$ of $TA$ over $TM$
\begin{equation}\label{liftptbaza} a_\alpha^C=c_\alpha,\;a_\alpha^V=
\frac{\partial}{\partial\xi^\alpha}.\end{equation} Notice also the
following properties
\begin{equation}\label{Cfunctie} (fa)^C=fa^C+f^Ca^V,\;(fa)^V=fa^V
\hspace{5mm}(f\in C^\infty(M),\,f^C=\dot{x}^i\frac{\partial
f}{\partial x^i}).\end{equation}

Now, assume that $(A,\sharp_A,[\,,\,]_A)$ is a Lie algebroid.
Define an anchor $\sharp_{TA}:TA\rightarrow TTM$ by putting
\begin{equation}\label{ancoraTA}
\sharp_{TA}c_\alpha=(\sharp_Aa_\alpha)^C,\;
\sharp_{TA}\frac{\partial}{\partial\xi^\alpha}=(\sharp_Aa_\alpha)^V,
\end{equation}
where in the right hand side the lifts are those of the Lie
algebroid $TM$. A simple calculation shows that
$\sharp_{TA}=\phi\circ(\sharp_A)_*$ with $\phi$ defined by
(\ref{flip}).

Furthermore, define a bracket of cross sections of $\pi_*$ by
putting
\begin{equation}\label{crosetTA}
[a_\alpha^V,a_\beta^V]_{TA}=0,\;[a_\alpha^C,a_\beta^V]_{TA}=
[a_\alpha,a_\beta]_A^V,\;[a_\alpha^C,a_\beta^C]_{TA}=[a_\alpha,a_\beta]_A^C,
\end{equation} which yields brackets of the elements of the basis of
cross sections of $TA\rightarrow TM$, and by extending
(\ref{crosetTA}) to arbitrary cross sections via the axioms of a Lie
algebroid. One can check that the results of (\ref{crosetTA}) hold
for the lifts of arbitrary cross sections of $A$ and this justifies
the independence of the bracket of the choice of the basis
$a_\alpha$. Then, the axioms of a Lie algebroid hold for the basic
cross sections $(c_\alpha,\partial/\partial\xi^\alpha)$, whence, the
axioms also hold for any cross sections.

Thus, we have obtained a well defined structure of a Lie algebroid
on the bundle $\pi_*:TA\rightarrow TM$. This Lie algebroid is
called the {\it tangent Lie algebroid} of $A$ \cite{MX}.

Now, we prove the announced result
\begin{prop}\label{tgfoliat} If the subalgebroid $B\subseteq A$ is
a foliation of $A$ over $(M,\mathcal{F})$, then the tangent Lie
algebroid $TB$ is a foliation of the Lie algebroid $TA$ over
$(TM,T\mathcal{F})$.\end{prop}
\begin{proof} We may use adapted local coordinates $(x^a,y^u)$ on $M$, such that
$ \mathcal{F}$ has the local equations $x^a=0$, and local bases
$(a_h)\equiv(b_\alpha,q_\kappa)$ of $\Gamma A$ such that
$(b_\alpha)$ is a local basis of $\Gamma B$ and $(q_\kappa)$ are
$B$-foliated cross sections of $A$. If we use formulas
(\ref{liftptbaza}), (\ref{ancoraTA}), (\ref{crosetTA}), we see
that the complete and vertical lifts of cross sections of
$B\subseteq A$ to $TA$ and to $TB$ coincide, and that the bracket
and anchor of $TB$ are the restrictions of the bracket and anchor
of $TA$ to $TB\subseteq TA$. Therefore, $TB$ is a Lie subalgebroid
of $TA$ over $TM$.

Furthermore, on $TM$ we have local coordinates
$(x^a,y^u,\dot{x}^a,\dot{y}^u)$ and (\ref{transTA}) shows that
$x^a=const.,\dot{x}^a=const.$ are the equations of a foliation
$T\mathcal{F}$ of the manifold $TM$ that consists of the tangent
vectors of the leaves of $\mathcal{F}$. Since
$\sharp_{TB}=\phi\circ(\sharp_B)_*$ where $\phi$ is the flip
diffeomorphism and $im\,\sharp_B=T\mathcal{F}$, we get
$im\,\sharp_{TB}=TT\mathcal{F}$.

Finally, on one hand, $(b_\alpha^C,b_\alpha^V)$ is a local basis
of cross sections of $TB$ and, on the other hand, using
(\ref{crosetTA}), it is easy to check that
$(q_\kappa^C,q_\kappa^V)$ are $TB$-foliated, local cross sections
of $TA\rightarrow TM$. Therefore, since
$(b_\alpha^C,q_\kappa^C,b_\alpha^V,q_\kappa^V)$ is a local basis
of cross sections of $TA\rightarrow TM$, the pair $(TA,TB)$
satisfies the foliated-generation condition.\end{proof}
\section{The dual Poisson structure}
It is well known that a Lie algebroid structure on the vector
bundle $\pi:E\rightarrow M$ is equivalent with a specific Poisson
structure on the total space of the dual bundle $E^*$, called the
{\it dual Poisson structure}, defined by the following brackets of
basic and fiber-linear functions:
\begin{equation}\label{PoissonE*} \{f_1\circ\pi,f_2\circ\pi\}=0,\;
\{f\circ\pi,l_s\}=-((\sharp_E
s)f)\circ\pi,\;\{l_{s_1},l_{s_2}\}=l_{[s_1,s_2]_E},
\end{equation} where $f,f_1,f_2\in C^\infty(M)$,
$s,s_1,s_2\in\Gamma E$ and $l_s$ is the evaluation of the fiber of
$E^*$ on $s$. If $x^a$ are local coordinates on $M$ and $\eta_h$
are the fiber-coordinates on $E^*$ with respect to the dual of a
local basis $(e_h)$ of cross sections of $E$, then the
corresponding Poisson bivector field is
\begin{equation}\label{Pi}
\Pi=\frac{1}{2}\alpha_{hk}^l\eta_l\frac{\partial}{\partial\eta_h}
\wedge\frac{\partial}{\partial\eta_k}
+\alpha_h^a\frac{\partial}{\partial\eta_h}\wedge\frac{\partial}{\partial
x^a},\end{equation} where the coefficients $\alpha$ are defined by
the expressions
\begin{equation}\label{coefPi}\sharp_E(e_h)=\alpha_h^a\frac{\partial}{\partial
x^a},\;[e_h,e_k]_E=\alpha_{hk}^le_l.\end{equation} Conversely,
formulas (\ref{Pi}), (\ref{coefPi}) produce an anchor $\sharp_E$
and a bracket $[\,,\,]_E$ and the Poisson condition $[\Pi,\Pi]=0$
implies the Lie algebroid axioms.

A Poisson structure of the form (\ref{Pi}) will be called a
fiber-linear structure, although this name is not totally
appropriate since it does not describe the form of the second term
of (\ref{Pi}). The characteristic property of a bivector field of
the form (\ref{Pi}) is that the Poisson bracket of two
fiber-polynomials of degrees $h,k$ is a fiber-polynomial of degree
$h+k-1$.

In this section we establish properties of the dual Poisson
structure of a foliated Lie algebroid.
\begin{prop}\label{dualfol1} Let $E$ be a foliated vector bundle
over $(M,\mathcal{F})$ with the dual bundle $E^*$, which has the
covering foliation $\mathcal{F}^{E^*}$ of $\mathcal{F}$. A
transversal-Lie algebroid structure on $E$ is equivalent with a
fiber-linear Poisson algebra structure on the space of
$\mathcal{F}^{E^*}$-foliated functions on $E^*$.\end{prop}
\begin{proof}
Let $(x^a,y^u)$ be $ \mathcal{F}$-adapted local coordinates on $M$
and $(e_h)$ be a local foliated basis of $E$. Then, the space of
foliated functions $C^\infty_{fol}(E^*,\mathcal{F}^{E^*})$ is
locally spanned by $(x^a,\eta_h)$ and we get the Poisson algebra
structure required by the proposition using formulas
(\ref{PoissonE*}) for $f,f_1,f_2\in
C_{fol}^\infty(M,\mathcal{F})$, $s,s_1,s_2\in\Gamma_{fol} E$.
Conversely, the Poisson structure (\ref{PoissonE*}) on
$C^\infty_{fol}(E^*,\mathcal{F}^{E^*})$ produces a transversal-Lie
algebroid on $E$ in the same way as in the case of a Lie
algebroid.\end{proof}
\begin{prop}\label{dualfol2} A foliated Lie algebroid $(A,B)$ over
$(M,\mathcal{F})$ is equivalent with a couple $(E^*,\Lambda)$,
where $E^*$ is a vector subbundle of $A^*$ endowed with a foliated
structure and the corresponding $ \mathcal{F}$-covering foliation
$\mathcal{F}^{E^*}$, and $\Lambda$ is a fiber-linear Poisson
structure on the manifold $A^*$ with the following properties:

1) $\sharp_\Lambda(ann TE^*)=T\mathcal{F}^{E^*}$,

2) $\Lambda|_{E^*}$ induces a well defined Poisson algebra structure
on $C^\infty_{fol}(E^*,\mathcal{F}^{E^*})$,

3) the total bundle manifold of $A^*/E^*$ has a Poisson structure
$P$ such that the projection $(A^*,\Lambda)\rightarrow (A^*/E^*,P)$
is a Poisson mapping.\end{prop}
\begin{proof} For the foliated Lie algebroid $(A,B)$, we
take $E^*=ann\,B$, which is foliated since its dual bundle is
$E=A/B$ and it is foliated. Furthermore, we take the dual Poisson
structure $\Lambda$ of the Lie algebroid $A$.

In order to write down the Poisson bivector field $\Lambda$ we
define convenient local coordinates on the manifold $A^*$ as
follows. We take $ \mathcal{F}$-adapted local coordinates
$(x^a,y^u)$ on $M$, we choose a splitting
\begin{equation}\label{ABC} A=B\oplus C,
\end{equation} and we take a local basis of $\Gamma A$ that
consists of the basis $b_h$ of $\Gamma B$ and the basis $a_q$ of
$\Gamma C$ where $a_q$ are foliated cross sections (see Lemma
\ref{bazefol}). For $A^*$, we have the dual bases $(b^{*h},a^{*q})$
and corresponding fiber coordinates $(\eta_h,\zeta_q)$.

Since $(A,B)$ is a foliated Lie algebroid, the anchor and bracket
have the following local expression
\begin{equation}\label{crosetbazeext}\begin{array}{l}
\sharp_A(b_h)=\beta_h^u\frac{\partial}{\partial y^u},\;
\sharp_A(a_q)=\alpha_q^a\frac{\partial}{\partial x^a}
+\alpha_q^u\frac{\partial}{\partial
y^u}\;(rank(\beta^u_h)=dim\,\mathcal{F}),\vspace{2mm}\\

[b_h,b_k]_A=\sum\beta_{hk}^lb_l,\,[b_h,a_q]_A=\gamma_{hq}^lb_l,\;
[a_p,a_q]_A=\alpha_{pq}^lb_l+\alpha_{pq}^sa_s,\end{array}\end{equation}
where $\alpha^a_q,\alpha_{pq}^s$ are foliated functions, i.e.,
locally, these are functions on the coordinates $x^a$.

Then, the general formula (\ref{Pi}) implies
\begin{equation}\label{PoissonpeA*} \Lambda=
\frac{1}{2}\beta_{hk}^l\eta_l\frac{\partial}{\partial\eta_h}
\wedge\frac{\partial}{\partial\eta_k}
+\gamma_{hq}^l\eta_l\frac{\partial}{\partial\eta_h}
\wedge\frac{\partial}{\partial\zeta_q}
+\frac{1}{2}(\alpha_{pq}^l\eta_l+\alpha_{pq}^s\zeta_s)
\frac{\partial}{\partial\zeta_p}
\wedge\frac{\partial}{\partial\zeta_q}\end{equation}
$$
+\alpha_q^a\frac{\partial}{\partial\zeta_q}\wedge\frac{\partial}{\partial
x^a}
+\beta_h^u\frac{\partial}{\partial\eta_h}\wedge\frac{\partial}{\partial
y^u}
+\alpha_q^u\frac{\partial}{\partial\zeta_q}\wedge\frac{\partial}{\partial
y^u}.$$

The total space of the subbundle $E^*\subseteq A^*$ has the local
equations $\eta_h=0$, whence
$$\sharp_\Lambda(ann\,TE^*)=span\{i(d\eta_h)\Lambda|_{\eta^h=0}\}
=span\{\beta^u_h\frac{\partial}{\partial
y^u}\}=T\mathcal{F}^{E^*},$$ which is property 1).

The space $C^\infty_{fol}(E^*,
\mathcal{F}^{E^*})$ is locally generated by $x^a,\zeta_q$
$(mod.\,\eta_h=0)$ and a function of local expression
$f(x^a,\zeta_q)$ extends to a function of local expression
$\tilde{f}(x^a,\zeta_q,\eta_h)$ in a neighborhood of $E^*$ in
$A^*$. Formula (\ref{PoissonpeA*}) restricted to $\eta_h=0$ shows
that, for $f_1,f_2\in C^\infty_{fol}(E^*,
\mathcal{F}^{E^*})$, $\Lambda|_{E^*}(d\tilde{f}_1,d\tilde{f}_2)$
depends only on $f_1,f_2$ and that
$$\{x^a,x^b\}_{\Lambda|_{E^*}}=0,\,
\{\zeta_q,x^a\}_{\Lambda|_{E^*}}=\alpha^a_q(x^b),\,
\{\zeta_p,\zeta_q\}_{\Lambda|_{E^*}}=\alpha^s_{pq}(x^b)\zeta_s.$$
Therefore, $\Lambda$ satisfies property 2).

Since $A^*/E^*\approx B^*$, we may use as local coordinates on the
manifold $A^*/E^*$ the coordinates $(x^a,y^u,\eta_h)$. Then,
\begin{equation}\label{eqP}P=\frac{1}{2}\beta_{hk}^l\eta_l
\frac{\partial}{\partial\eta_h}
\wedge\frac{\partial}{\partial\eta_k} +\beta_h^u
\frac{\partial}{\partial\eta_h}\wedge\frac{\partial}{\partial
y^u}\end{equation} is a bivector field on $A^*/E^*$, which is
Poisson because it has the same expression as the dual Poisson
structure of the Lie algebroid $B$. Property 3) is obviously
satisfied.

Conversely, let $A$ be a vector bundle over $(M,\mathcal{F})$ such
that there exists a foliated subbundle $E^*$ of the dual bundle
$A^*$ and a fiber-linear Poisson bivector field $\Lambda$ of $A^*$
with the properties 1), 2), 3). Define the subbundle $B=ann
E^*\subseteq A$ and take a decomposition $A=B\oplus C$ and local
bases of cross sections $b_h\in B,a_q\in C$. Moreover, $C\approx
A/B\approx E=E^{**}$ is an $
\mathcal{F}$-foliated bundle and $a_q$ may be assumed to be $
\mathcal{F}$-foliated cross sections. The corresponding local
coordinates $(x^a,y^u,\eta_h,\zeta_q)$ on $A^*$ are similar to
those used in (\ref{PoissonpeA*}) and, a priori, $\Lambda$ is of
the form $$ \Lambda=
\frac{1}{2}(\beta_{hk}^l\eta_l+ \beta_{hk}^s\zeta_s)
\frac{\partial}{\partial\eta_h}
\wedge\frac{\partial}{\partial\eta_k}
+(\gamma_{hq}^l\eta_l+\gamma_{hq}^s\zeta_s)\frac{\partial}{\partial\eta_h}
\wedge\frac{\partial}{\partial\zeta_q}$$
$$+\frac{1}{2}(\alpha_{pq}^l\eta_l+\alpha_{pq}^s\zeta_s)
\frac{\partial}{\partial\zeta_p}
\wedge\frac{\partial}{\partial\zeta_q}
+\alpha_q^a\frac{\partial}{\partial\zeta_q}\wedge\frac{\partial}{\partial
x^a}
+\beta_h^u\frac{\partial}{\partial\eta_h}\wedge\frac{\partial}{\partial
y^u}$$
$$+\alpha_q^u\frac{\partial}{\partial\zeta_q}\wedge\frac{\partial}{\partial
y^u}+\lambda_h^a\frac{\partial}{\partial\eta_h}\wedge\frac{\partial}{\partial
x^a}.$$ If we put $\eta_h=0$ and compute $i(d\eta_h)\Lambda$ we
see that property 1) implies
$$\beta^s_{hk}=0,\,\gamma^s_{hq}=0,\,\lambda^a_h=0,\,rank(\beta^u_h)=dim\,
\mathcal{F}.$$ Furthermore, property 2) implies local expressions
$\alpha^a_q=\alpha^a_q(x^b),\,\alpha^s_{pq}=\alpha^s_{pq}(x^b)$.
Finally, the projection ($A^*\rightarrow A^*/E^*$)-related
bivector field $P$ of property 3) is uniquely defined by
$\Lambda$; it must be given by (\ref{eqP}) and it satisfies the
Poisson condition. Therefore, $\Lambda$ reduces to the form
(\ref{PoissonpeA*}) and the corresponding Lie algebroid structure
of $A$ is given by formulas (\ref{crosetbazeext}). Accordingly,
the pair $(A,B)$ is a foliated Lie algebroid. (The
foliated-generation condition is implied by the fact that
$[b_h,a_q]\in B$.)\end{proof}

A. Vaintrob \cite{Vb} gave an interpretation of the dual Poisson
structure (\ref{Pi}) as an odd, homological super-vector field,
which led to important results on Lie algebroids seen as
homological super-vector fields. We shall indicate the properties
of this field in the foliated case.

The {\it parity-changed vector bundle} of $A^*$ (in our case $A^*$
is the dual of a Lie algebroid, but, the definition applies to any
vector bundle) is a supermanifold $\Pi A^*$, with local even
coordinates $x^i$ and local odd coordinates $\bar{\eta}^\alpha$
associated to local, $A^*$-trivializing, coordinate neighborhoods
$U$ on $M$, with the following local transition functions:
$$\tilde{x}^i=\tilde{x}^i(x^j),\,\tilde{\bar{\eta}}^\alpha =
b^\alpha_\beta\bar{\eta}^\beta,$$ where the first equations are
the same as the transition functions of the local coordinates on
$M$ and the matrix $(b^\alpha_\beta)$ is defined by a change of
local bases of cross sections of $A$,
$a_\beta=b^\alpha_\beta\tilde{a}_\alpha$. The supertangent bundle
$T\Pi A^*$ has the local bases $({\partial}/{\partial x^i},\,
{\partial}/{\partial
\bar{\eta}^\alpha})$ with the local transition functions
$$\frac{\partial}{\partial \tilde{x}^i}=
\frac{\partial x^j}{\partial \tilde{x}^i}\frac{\partial}{\partial
x^j},\, \frac{\partial}{\partial \bar{\eta}^\beta}=
b^\alpha_\beta\frac{\partial}{\partial\tilde{\bar{\eta}}^\alpha}.$$

These transition functions show that the correspondence
$$x^i\mapsto x^i,\,\frac{\partial}{\partial x^i}\mapsto
\frac{\partial}{\partial x^i},\,\eta_\alpha\mapsto
\frac{\partial}{\partial \bar{\eta}^\alpha},\,
\frac{\partial}{\partial {\eta}_\alpha}\mapsto\bar{\eta}^\alpha$$
produces a well defined mapping of the fiber-wise polynomial
multivector fields on the manifold $A^*$ to super-multivector
fields on $\Pi A^*$. In particular, the bivector field $\Pi$
defined by (\ref{Pi}) is send to Vaintrob's super-vector field
\begin{equation}\label{Vfield} V=\frac{1}{2}\alpha_{hk}^l\bar{\eta}^h\bar{\eta}^k
\frac{\partial}{\partial\bar{\eta}^l}
+\alpha_h^a\bar{\eta}^h\frac{\partial}{\partial x^a}\end{equation}
and the Poisson condition $[\Pi,\Pi]=0$ becomes the homology
condition $[V,V]=0$. To put terminology in agreement with the
theory of Lie algebroids, the super-vector field $V$ will be
called {\it transitive} if the $M$-{\it rank} $\rho$ of $V$
defined by $\rho=rank(\alpha_h^a)$ is equal to the dimension of
$M$.

Furthermore, if $M$ is endowed with the foliation $\mathcal{F}$
and $A$ has a subbundle $B$, we have coordinates $(x^a,y^u)$ and
bases $(b_h,a_q)$ like in the proof of Proposition \ref{dualfol2},
which yield coordinates $(x^a,y^u,\bar{\eta}^h,\bar{\zeta}^q)$ of
$\Pi A^*$ with transition functions of the form
\begin{equation}\label{supertransfol}\tilde{x}^a=\tilde{x}^a(x^b),\,
\tilde{y}^u=\tilde{y}^u(x^a,y^v),\,  \tilde{\bar{\eta}}^h =
b^h_k\bar{\eta}^k+b^h_q\bar{\zeta}^q,\,\tilde{\bar{\zeta}}^q=b^q_p\bar{\zeta}^p.
\end{equation} where the various matrices $b$ come from a change
of the local basis. The transition functions (\ref{supertransfol})
show that $\Pi A^*$ is endowed with the superfoliation
$$\bar{\mathcal{F}}=span\left\{\frac{\partial}{\partial y^u},
\frac{\partial}{\partial\bar{\eta}^h}\right\}.$$ Then, if we
consider the vector bundle $E= A/B$, we have $E^*=ann\,B\subseteq
A^*$ and $\Pi(A^*/E^*)$ may be seen as the submanifold of $\Pi
A^*$ that has the local equations $\bar{\zeta}^q=0$. On the other
hand, $\Pi E^*$ is a supermanifold with the local coordinates
$(x^a,y^u,\bar{\zeta}^q)$ and there exists a natural projection
$\psi:\Pi A^*\rightarrow\Pi E^*$ locally defined by
$$(x^a,y^u,\bar{\zeta}^q,\bar{\eta}^h)\mapsto
(x^a,y^u,\bar{\zeta}^q).$$

In Vaintrob's terminology a super-vector field on $\Pi A^*$ that
corresponds to a bivector field on $A^*$ is said to be of degree
$1$. In our case, a super-vector field of degree $1$ has the form
\begin{equation}\label{superW}
W=\frac{1}{2}\beta_{hk}^l\bar{\eta}^h\bar{\eta}^k
\frac{\partial}{\partial\bar{\eta}^l} +
\frac{1}{2}\beta_{hk}^s\bar{\eta}^h\bar{\eta}^k
\frac{\partial}{\partial\bar{\zeta}^s}
+ \gamma_{hq}^l\bar{\eta}^h\bar{\zeta}^q
\frac{\partial}{\partial\bar{\eta}^l}\end{equation}
$$+ \gamma_{hq}^s\bar{\eta}^h\bar{\zeta}^q
\frac{\partial}{\partial\bar{\zeta}^s}
+\frac{1}{2}\alpha_{pq}^l\bar{\zeta}^p\bar{\zeta}^q
\frac{\partial}{\partial\bar{\eta}^l}
+\frac{1}{2}\alpha_{pq}^s\bar{\zeta}^p\bar{\zeta}^q
\frac{\partial}{\partial\bar{\zeta}^s}$$
$$+\alpha_q^a\bar{\zeta}^q\frac{\partial}{\partial x^a}
+\beta_h^u\bar{\eta}^h\frac{\partial}{\partial y^u}
+\alpha_q^u\bar{\zeta}^q\frac{\partial}{\partial y^u}
+\lambda_h^a\bar{\eta}^h\frac{\partial}{\partial x^a}.$$
\begin{prop}\label{algfolcuV} A foliated Lie algebroid
$(A,B)$ over $(M,\mathcal{F})$ is equivalent with a couple
$(E^*,W)$, where $E^*$ is an $\mathcal{F}$-foliated vector
subbundle of $A^*$ and $W$ is a homological super-vector field of
degree $1$ on $\Pi A^*$ such that

i) $W$ is foliated (projectable) with respect to the super-foliation
$\bar{\mathcal{F}}$ of $\Pi A^*$ associated to the pair
$(\mathcal{F},B)$,

ii) $W|_{\Pi(A^*/E^*)}\in T\bar{\mathcal{F}}\cap T\Pi(A^*/E^*)$ and
is a homological super-vector field of degree $1$ on $\Pi(A^*/E^*)$,
transitive over the leaves of $\mathcal{F}$.\end{prop}
\begin{proof}
Formula (\ref{PoissonpeA*}) shows that the Vaintrob super-vector
field $W$ of a foliated Lie algebroid $(A,B)$ is a homological
super-vector field (\ref{superW}) that satisfies the conditions
\begin{equation}\label{eqaux} \begin{array}{l}
\gamma_{hq}^s=0,\,\alpha_{pq}^s=\alpha_{pq}^s(x^b),\,
\alpha_q^a=\alpha_q^a(x^b),\vspace{2mm}\\ \lambda_h^a=0,\,\beta_{hk}^s=0,
\,rank(\beta_h^u)=dim\,\mathcal{F},\end{array}\end{equation}
and where the super-vector field
$$W|_{\Pi(A^*/E^*)}=\frac{1}{2}\beta_{hk}^l\bar{\eta}^h\bar{\eta}^k
\frac{\partial}{\partial\bar{\eta}^l}
+\beta_h^u\bar{\eta}^h\frac{\partial}{\partial y^u}\in
T\bar{\mathcal{F}}\cap T\Pi(A^*/E^*)$$ is a homological
super-vector field on $\Pi(A^*/E^*)$. The condition
$rank(\beta_h^u)=dim\,\mathcal{F}$ is what we meant by the
transitivity of $W$ over the leaves of $
\mathcal{F}$. Conversely, if we start with the data $(A,E^*,W)$,
$A$, which has the Lie algebroid structure of Vaintrob
super-vector field $W$, has the subbundle $B=ann\,E^*$ and $\Pi
A^*$ has the superfoliation $\bar{\mathcal{F}}$. Then, $W$ may be
represented under the form (\ref{superW}). Condition i), which
means that the terms of (\ref{superW}) that contain
$(\partial/\partial x^a,\partial/\partial\bar{\zeta}^q)$ depend on
the coordinates $(x^a,\bar{\zeta}^q)$ only, implies the first four
conditions (\ref{eqaux}) and condition ii) ensures the fact that
$B$, with the Lie algebroid structure of Vaintrob super-vector
field $W|_{\Pi(A^*/E^*)}$, is a foliation of $A$.\end{proof}
\section{Cohomology and deformations}
The cohomology $H^*(A)$ of a Lie algebroid $A$ is that of the
differential graded algebra $(\Omega(A)=\Gamma\wedge A^*,d_A)$ of
$A$-{\it forms} where, if $\lambda\in\Omega^k(A)$, then
\begin{equation}\label{difext} d_A\lambda(s_1,...,s_{k+1})=
\sum_{j=1}^{k+1}(-1)^{j+1}\sharp_As_j(\lambda(s_1,...,s_{j-1},s_{j+1},...,s_{k+1}))
\end{equation}
$$+\sum_{j<l=1}^{k+1}(-1)^{j+l}\lambda([s_j,s_l]_A,s_1,...,s_{j-1},s_{j+1},...,
s_{l-1},s_{l+1},...,s_{k+1}).$$ The anchor of the Lie algebroid
produces a homomorphism of differential graded algebras
$\sharp'_A:(\Omega(M),d)\rightarrow(\Omega(A),d_A)$ defined by
$$(\sharp'\theta)(s_1,...,s_{k})=
\theta(\sharp s_1,...,\sharp s_{k}),\;\;\theta\in\Omega^k(M)$$
with the corresponding cohomology homomorphism $\sharp^*_A:H^*_{de
Rham}(M)\rightarrow H^*(A)$.

The general pattern of the cohomology of a foliated Lie algebroid
$(A,B)$ over a foliated base manifold $(M,\mathcal{F})$ is the
same like that of a foliated manifold
\cite{V1973}, which is why we only sketch it briefly here. Using
the decomposition $A=B\oplus C$ and the bases $(b_h,a_l)$,
$(b^{*h},a^{*l})$ where $b_h\in\Gamma B,a_l\in\Gamma C$ and $a_l$
are $B$-foliated, we get a bigrading of the $A$-forms such that a
form of {\it type} (bi-degree) $(s,r)$ ($s$ is the $C$-grade and
$r$ is the $B$-grade) has the local expression
\begin{equation}\label{biexpresie}
\lambda=\frac{1}{s!r!}\lambda_{l_1...l_sh_1...h_r}(x^a,y^u)a^{*l_1}
\wedge...\wedge a^{*l_s}\wedge b^{*h_1}\wedge...\wedge
b^{*h_r}.\end{equation} We will denote by $\Omega^{s,r}(A)$ the
space of $A$-forms of type $(s,r)$.

Since $B$ is bracket closed, one has $d_Aa^{*l}(b_h,b_k)=0$ and we
see that $d_A\lambda$ has only components of the type
$(s+1,r),(s,r+1),(s+2,r-1)$. Hence there exists a decomposition
\begin{equation}\label{descdA}
d_A=d'_A+d''_A+\partial_A,\end{equation} where the terms are
differential operators of types $(1,0),(0,1),(2,-1)$,
respectively. The property $d_A^2=0$ is equivalent with
\begin{equation}\label{propdA} \begin{array}{c}(d''_A)^2=0, d'_Ad''_A +
d''_Ad'_A=0,\;(\partial_A)^2=0,\vspace{2mm}\\ d'_A\partial_A +
\partial_Ad'_A=0,\;(d'_A)^2 + d''_A\partial_A +
\partial_Ad''_A=0. \end{array}\end{equation}

The cohomology of $(\Omega^{0\,\centerdot}(A),d''_A)$ is the
cohomology of the Lie algebroid $B$. On the other hand, an
$A$-form $\lambda$ is said to be $B$-{\it foliated} if it is of
type $(s,0)$ and its local components $\lambda_{l_1,...,l_s}$ are
$\mathcal{F}$-foliated functions (this is an invariant condition
because the bases $(a_l)$ are $B$-foliated and have
$\mathcal{F}$-foliated transition functions). The operator $d'_A$
induces a cohomology of the $B$-foliated forms called the {\it
basic cohomology}. Formulas (\ref{propdA}) show that
$(\Omega(A),d_A)$ is a semi-positive double cochain complex
\cite{V1973}, hence, it produces a spectral sequence that connects
among the previously mentioned cohomologies. Another component of
the pattern is the truncated cohomology discussed in the Appendix
of \cite{V3}.
\begin{example}\label{coh1dif} {\rm In Example \ref{subreg} we
showed that for any Lie algebroid $A$ the pair $(J^1A,Hom(TM,A))$
is a foliated Lie algebroid over the foliation of the basis $M$ by
points, therefore, its cohomology has the pattern described above.
On the other hand, the cohomology of $J^1A$ may be interpreted as
the $1$-{\it differentiable cohomology} of $A$ in the sense of
Lichnerowicz, who studied the $1$-differentiable cohomology of
many infinite dimensional Lie algebras on manifolds (e.g., see
\cite{{L1},{L2}} and the references therein; the case studied in
\cite{L2} is that of the cotangent Lie algebroid of a Poisson manifold).
This cohomology is defined by $k$-cochains $\lambda\in
Hom_{\mathds{R}} (\Gamma(\wedge^kA),C^\infty(M))$, which are
differential operators of order $1$ in each argument, and by the
differential $d$ defined by formula (\ref{difext}). Since a
differential operator of order $1$ is the composition of a
$C^\infty(M)$-linear morphism by the first jet mapping $j^1$, the
$1$-differentiable cochains may be identified with usual cochains
of the Lie algebroid $J^1A$. Furthermore, since
$[j^1a_1,j^1a_2]_{J^1A}=j^1[a_1,a_2]_A$ $(a_1,a_2\in\Gamma A)$
\cite{CF}, $d$ for $1$-differentiable cochains may be identified
with $d_{J^1A}$. Hence, $H^k_{1-diff}(A)\approx H^k(J^1A)$ as
claimed. A complementary subbundle $C$ of the foliation
$Hom(TM,A)$ is given by the image of a splitting
$\sigma:A\rightarrow J^1A$ of the exact sequence (\ref{sirjeturi})
and the type decomposition (\ref{biexpresie}) corresponds to the
type decomposition used by Lichnerowicz.}\end{example}

Like in the case of a foliation \cite{V4}, we get
\begin{prop}\label{Poincare} {\rm (The $d''_A$-Poincar\'e Lemma)}
Let $(A,B)$ be a minimally foliated Lie algebroid. For any local
$A$-form $\lambda$ of type $(u,v)$ $(v>0)$ such that
$d''_A\lambda=0$, one has $\lambda=d''_A\mu$ where $\mu$ is a
local form of type $(u,v-1)$.\end{prop}
\begin{proof} We may proceed by
induction on $k$, where $k$ is the number defined by the condition
that the local expression of $\lambda$ does not contain the
elements $a^{*k+1},...,a^{*q}$ $(q=rank\,A-rank\,B)$ of the basis
used in (\ref{biexpresie}).

Since we are in the case of a minimally foliated algebroid,
$\sharp_B:B\rightarrow T\mathcal{F}$ is an isomorphism, which
induces an isomorphism between the graded differential algebra
$(\Omega^{0\,\centerdot}(A),d''_A)$ and the graded differential
algebra of the $\mathcal{F}$-leaf-wise forms on $(M,\mathcal{F})$.
Since the $d''$-Poincar\'e lemma holds for the latter
\cite{{V4},{V1973}}, it holds for the former as well, which
justifies the conclusion for $k=0$.

Now, we assume that the result holds for any $k<h-1$ $(h=1,...,q)$
and take
$$\lambda=a^{*h}\wedge\mu+\theta$$ where $\mu,\theta$ do not
contain $a^{*h},...,a^{*q}$. Since $d''_Aa^{*h}=0$, $d''_A\lambda=0$
if and only if $d''_A\mu=d''_A\theta=0$. Then, the induction
hypothesis yields the local expressions
$\mu=d''_A\tau,\theta=d''_A\sigma$, therefore,
$$\lambda=d''_A(-a^{*h}\wedge\tau+\sigma),$$ which means that the
result also holds for $k=h-1$ and we are done.
\end{proof}

Obviously, the operator $d''_A$ makes sense for $V$-valued
$A$-forms, where $V$ is any $ \mathcal{F}$-foliated vector bundle
on M, and the $d''$-Poincar\'e lemma (Proposition \ref{Poincare})
holds for such forms as well.

The following corollary is a standard result.
\begin{corol}\label{corolRham} Let $(A,B)$ be a minimally
foliated Lie algebroid and $V$ a foliated vector bundle over
$(M,\mathcal{F})$. Let $\Phi^k$ be the sheaf of germs of
$B$-foliated, $V$-valued $A$-forms. Then one has
$$ H^l(M,\Phi^k)\approx\frac{ker (d''_A:\Omega^{k,l}(A)\otimes V
\rightarrow\Omega^{k,l+1}(A)\otimes
V)}{im(d''_A:\Omega^{k,l-1}(A)\otimes V
\rightarrow\Omega^{k,l}(A)\otimes V)}.$$\end{corol}
\begin{rem}\label{Prelativ}
{\rm In the terminology of \cite{V3}, Proposition \ref{Poincare}
is a relative Poincar\'e lemma. Therefore the sheaf-interpretation
of the truncated cohomology given in
\cite{V3} holds for minimally foliated Lie algebroids.}\end{rem}

There exists a general construction of characteristic
$A$-cohomology classes of a Lie algebroid \cite{F}, which mimics
the Chern-Weil theory and goes as follows. Take a vector bundle
$V\rightarrow M$ of rank $r$. An $A$-connection (covariant
derivative) on $V$ is an $\mathds{R}$-bilinear operator
$\nabla:\Gamma A\times
\Gamma V\rightarrow\Gamma V$ that is $C^\infty(M)$-linear in the
first argument and satisfies the condition
$$\nabla_a(fv)=f\nabla_av+((\sharp_Aa)f)v,\hspace{3mm}(a\in\Gamma
A,v\in\Gamma V,f\in C^\infty(M)).$$ The curvature operator of
$\nabla$ is $$R_\nabla(a_1,a_2)=\nabla_{a_1}\nabla_{a_2}
-\nabla_{a_2}\nabla_{a_1}-\nabla_{[a_1,a_2]_A}.$$ $R_\nabla$ is
$C^\infty(M)$-bilinear and may be seen as a $V\otimes V^*$-valued
$A$-form of degree $2$. Then, $\forall\phi\in
I^k(Gl(r,\mathds{R}))$, the space of real, ad-invariant,
symmetric, $k$-multilinear functions on the Lie algebra
$gl(r,\mathds{R})$, the $A$-cohomology classes
$[\phi(R_\nabla)]=[\phi(R_\nabla,...,R_\nabla)]\in H^{2k}(A)$ are
the $A$-{\it principal characteristic classes} of $V$. The
characteristic classes are spanned by the classes
$[c_k(R_\nabla)]$, where $c_k$ is the sum of the principal minors
of order $k$ in $det(R_\nabla-\lambda Id)$, and do not depend on
the choice of the connection $\nabla$. Indeed, if
$\nabla^0,\nabla^1$ are two $A$-connections on $V$, Bott's proof
\cite{{Bott},{F}} gives
\begin{equation}\label{eqB1} \phi(R_{\nabla_1})-\phi(R_{\nabla^0})
=d_A\Delta(\nabla^0,\nabla^1)\phi, \end{equation} where
\begin{equation}\label{eqB2} \Delta(\nabla^0,\nabla^1)\phi=
k\int_0^1\phi(\nabla^1-\nabla^0,\underbrace{R_{(t\nabla^1+(1-t)\nabla^0)},
...,R_{(t\nabla^1+(1-t)\nabla^0)}}_{(k-1)-{\rm
times}})dt,\end{equation} hence the cohomology classes defined by
$\phi(R_\nabla^0),\phi(R_\nabla^1)$ coincide. In particular, the
use of an orthogonal $A$-connection (such that $\nabla g=0$ for a
Euclidean metric $g$ on $A$) gives $[c_{2h-1}(R_\nabla)]=0$. The
classes $[c_{2h}(R_\nabla)]$ are the $A$-{\it Pontrjagin classes}
of $V$. For $V=A$ one gets the principal characteristic classes of
the Lie algebroid $A$.

In the case of a foliated Lie algebroid $(A,B)$ over
$(M,\mathcal{F})$ one can also mimic the construction of secondary
characteristic classes. We define a {\it Bott $A$-connection} on a
foliated vector bundle $V$ by asking it to satisfy the condition
$\nabla_bv=0$, $\forall b\in B$, $\forall v\in\Gamma_{fol}V$. Bott
$A$-connections always exist: take an arbitrary $A$-connection
$\nabla'$ on $V$ and a decomposition $A=B\oplus C$, then, define
$$\nabla_{b+c}(fv)=((\sharp_Ab)f)v+\nabla'_c(fv)\hspace{2mm}
(f\in C^\infty(M),v\in\Gamma_{fol}V).$$ The curvature of a Bott
$A$-connection satisfies the property
\begin{equation}\label{Rbasic} R_\nabla(b_1,b_2)v=0, \forall b_1,b_2\in B,v\in
V.\end{equation} Indeed, due to $C^\infty(M)$-linearity with
respect to $v$, it suffices to check the property for
$v\in\Gamma_{fol}V$, which is trivial.
\begin{prop}\label{0phenomen} {\rm(The Bott Vanishing
Phenomenon)} If $k>q=rank\,A-rank\,B$, $\forall\phi\in
I^k(Gl(r,\mathds{R}))$, the characteristic class defined by $\phi$
vanishes. Equivalently, $Pont_A^h(V)=0$ for $h>2q$, where
$Pont_A^h(V)$ is the set of elements of total cohomological degree
$h$ generated by the $A$-Pontrjagin classes in the cohomology
algebra $H^*(A)$.\end{prop}
\begin{proof} Property (\ref{Rbasic}) shows that $R_\nabla$ belongs to
the ideal generated by $ann\,B$, hence, if $k>q$, the $A$-form
obtained by the evaluation of $\phi$ on $R_\nabla$ is
zero.\end{proof}
\begin{corol}\label{corolBott} Let $B$ be a Lie subalgebroid of the
Lie algebroid $A$. If $Pont_A^h(A/B)\neq0$ for some
$h>2(rank\,A-rank\,B)$, $B$ is not a foliation of $A$.\end{corol}

Now, if we take $\phi\in I^{2h-1}(Gl(r,\mathds{R}))$ with
$2h-1>q$, formula (\ref{eqB1}) yields
$d_A(\Delta(\nabla^0,\nabla^1)\phi)=0$ whenever $\nabla^0$ is an
orthogonal $A$-connection and $\nabla^1$ is a Bott $A$-connection
on the $\mathcal{F}$-foliated vector bundle $V$. The
$A$-cohomology classes $[\Delta(\nabla^0,\nabla^1)\phi]\in
H^{4h-3}(A)$ are {\it $A$-secondary characteristic classes of
$V$}. Like in \cite{L} one can prove that these classes do not
depend on the choice of $\nabla^0$ and remain unchanged if
$\nabla^1$ is subject to a deformation via Bott connections. In
the case $V=A/B$, the construction provides secondary
$A$-characteristic classes of the foliated Lie algebroid $(A,B)$.
\begin{prop}\label{imagini}
The $A$-secondary characteristic classes of $V$ are anchor images
of de Rham-secondary characteristic classes.\end{prop}
\begin{proof} We may define the $A$-secondary characteristic
classes by means of connections $\nabla^0,\nabla^1$ given by
$\nabla^0_a=\nabla^{'0}_{\sharp_Aa},\nabla^1_a=\nabla^{'1}_{\sharp_Aa}$,
where $\nabla'_0$ is a usual orthogonal connection on
$V\rightarrow M$ and $\nabla^{'1}$ is a Bott connection with
respect to the foliated Lie algebroid $(TM,T\mathcal{F})$. Then,
the conclusion holds at the level of forms.\end{proof}

Deformations of foliations were studied by many authors
\cite{{He},{KSp}}, etc., and one of the main results is that
deformations produce infinitesimal deformations, which are
cohomology classes of degree $1$. Here, we extend this result to
deformations of a foliation $B$ of a Lie algebroid $A$ over
$(M,\mathcal{F})$.
\begin{defin}\label{defdeform} {\rm
A {\it deformation} of the foliation $B$ of $A$ is a
differentiable family of foliations $B_t\subseteq A$ $(0\leq
t\leq1)$ where $B_0=B$, each $B_t$ is regular and
$\mathcal{F}_t=\sharp_A(B_t)$ is a deformation of the foliation
$\mathcal{F}=\mathcal{F}_0$.}\end{defin}

By differentiability of $B_t$ we mean that each $x\in M$ has a
neighborhood $U$ endowed with local cross sections
$(b_h(x,t))_{h=1}^p$ of $A$ that are differentiable with respect
to $(x,t)$ and define a basis of the local cross sections of the
subbundle $B_t$ for all $t\in[0,1]$. In particular, all the
subbundles $B_t$ have the same rank $p$.

The transition functions of the bases $(b_h(x,t))$ are of the form
\begin{equation}\label{transindeform}
\tilde{b}_h(x,t)=\beta_h^k(x,t)b_k(x,t).\end{equation} If we apply
the operator $(d/dt)_{t=0}$, we see that the correspondence
$$\eta^hb_h\mapsto \left.[\eta^h\frac{\partial b_h(x,t)}{\partial
t}\right|_{t=0}]_{{\rm mod.}\,B} \;\;(b_h=b_h(x,0))$$ provides a
well defined $E=A/B$-valued $1$-$B$-form $\Xi$.  The form $\Xi$
will be called the {\it infinitesimal deformation form} and it has
the following invariant expression
\begin{equation}\label{Xiglobal}
\Xi(b)=\left.[\frac{\partial\tilde{b}(x,t)}{\partial
t}\right|_{t=0}]_{{\rm mod.}\,B}
\end{equation} where $\tilde{b}(x,t)\in\Gamma B_t$ is an extension of $b\in B_x$;
the result does not depend on the choice of the extension.

In what follows we identify the cochain complexes
$(\Omega(B),d_B)$ and $(\Omega^{0\centerdot}(A),d''_A)$ by using a
complementary subbundle of $B$ in $A$; in particular, $\Xi$ is
also seen as a $(0,1)$-form.
\begin{prop}\label{inchidere} The infinitesimal deformation form $\Xi$ is
$d''_A$-closed.\end{prop}
\begin{proof} We begin by deriving a new
expression of $\Xi$. Let us fix a subbundle $C\subseteq A$ such
that $A=B_t\oplus C$ holds $\forall t\in[0,\epsilon)$ (where
$\epsilon$ is small enough) and consider a local basis $(b_h)$ of
$B$ and a local basis $(c_s)$ of $C$ consisting of $B$-foliated
cross sections; this yields a foliated basis $e_s=[c_s]_{{\rm
mod.}\,B}$ of $A/B$. Then, we have expressions
\begin{equation}\label{exprcuC}
b_h(x,t)=\lambda_h^k(x,t)b_k+\mu_h^s(x,t)c_s,
\hspace{2mm}rank(\lambda_h^k)=p,\end{equation} therefore,
$B_t$ also has a basis of the form
\begin{equation}\label{exprcuC2}
b'_h(x,t)=b_h+\theta_h^s(x,t)c_s,\hspace{3mm}\theta_h^s(x,0)=0.\end{equation}

The dual cobases of the bases $(b'_h(x,t),c_s)$ are given by
\begin{equation}\label{cobazet}
b^{'*h}(x,t)=b^{*h}(x,0),\;c^{*s}(x,t)=c^{*s}-\theta_h^s(x,t)b^{*h}
\end{equation} (check the duality conditions
$$<b^{'*h}(x,t),b'_k(x,t)>=\delta_k^h,\, <c^{*s}(x,t),b'_k(x,t)>=0,$$
$$<b^{'*h}(x,t),c_s>=0,\,<c^{*s}(x,t),c_u(x,t)>=\delta_u^s).$$

From (\ref{exprcuC2}) and (\ref{cobazet}) we get
$$\left.\left[\frac{\partial b'_h}{\partial
t}\right|_{t=0}\right]_{{\rm mod.}\,B}=\left.\frac{\partial
\theta^s_h}{\partial t}\right|_{t=0}e_s$$ and
\begin{equation}\label{exprXicuC} \Xi= \left.\frac{\partial \theta^s_h}{\partial
t}\right|_{t=0}b^{*h}\otimes e_s =\left.-\frac{\partial
c^{*s}(x,t)}{\partial t}\right|_{t=0}\otimes e_s.\end{equation}

If we denote by $\Xi^C$ the image of $\Xi$ by the vector bundle
isomorphism $E\approx C$, which following (\ref{exprXicuC}) means
that
\begin{equation}\label{XiC}\Xi^C=\left.-\frac{\partial c^{*s}(x,t)}{\partial
t}\right|_{t=0}\otimes c_s,\end{equation} the conclusion of the
proposition is equivalent with
$$d_A\Xi^C(b_h,b_k)=0,$$ which we see to hold by means of the
following calculation. From (\ref{XiC}) and using the definition
of $d_A$ we have
$$d_A\Xi^C(b_h,b_k)=-\left.\frac{\partial
[d_Ac^{*s}(x,t)(b'_h(x,t),b'_k(x,t))]}{\partial
t}\right|_{t=0}c_s$$ $$=\left.\frac{\partial
[c^{*s}(x,t)([b'_h(x,t),b'_k(x,t)]_A)]}{\partial
t}\right|_{t=0}c_s=\left.\frac{\partial}{\partial
t}\right|_{t=0}(pr^t_C[b'_h(x,t),b'_k(x,t)]_A),$$ where $pr^t_C$
is the projection defined by the decomposition $A=B_t\oplus C$.
Since $B_t$ is closed by $A$-brackets, the final term of the
previous calculation is equal to zero.\end{proof}

Continuing to follow the analogy with foliation theory, we give
the following definition.
\begin{defin}\label{defdeformtr} {\rm A {\it trivial deformation}
of the foliation $B$ of the Lie algebroid $A$ is a deformation
$B_t$ such that $B_t=\Psi_t(B)$, where $(\Psi_t,\psi_t)$
$(\psi_t:M\rightarrow M)$ is a family of automorphisms of $A$ of
the form $\Psi_t=exp(ta),\psi_t=exp(t\sharp_Aa)$, $a\in\Gamma A$.}
\end{defin}

We recall that $\Psi_t=exp(ta)$ is a (local) $1$-parameter group
with respect to addition on $t$, which is characterized by the
property \cite{M}
\begin{equation}\label{expa}\left.
\frac{d}{dt}\right|_{t=0}(exp(ta))a'=[a,a']_A, \hspace{3mm} \forall
a,a'\in\Gamma A.\end{equation}
\begin{prop}\label{Xitrivial} The infinitesimal deformation form
$\Xi$ of a trivial deformation is $d''_A$-exact.\end{prop}
\begin{proof} For a trivial deformation we may use the local bases
$$b_h(x,t)=exp(ta)(b_h(exp(-t\sharp_Aa)(x),0),$$ which implies
$$\left.\frac{\partial b_h(x,t)}{\partial t}\right|_{t=0}=[a,b_h]_A(x).$$
Now, using a decomposition $a=\alpha^hb_h+\beta^{s}c_s$, where
$(b_h,c_s)$ is the basis used in (\ref{exprcuC}), the original
definition of $\Xi$ shows that the corresponding form $\Xi^C$ is
given by
\begin{equation}\label{auxeq}
\Xi^C(\eta^hb_h)=pr_C(\eta^h[a,b_h]_A)=-\eta^h(\sharp_Ab_h(\beta^{s}))c_s
=-d''_A(pr_Ca),\end{equation} which is equivalent with the
required result.\end{proof}
\begin{defin}\label{deformechiv} {\rm Let $(A,B)$ be a foliated
Lie algebroid. Two deformations $B_t,\bar{B}_t$ of $B$ are {\it
equivalent} if there exists a cross section $a\in\Gamma A$ such
that $\bar{B}_t=exp(ta)(B_t)$, i.e., if they can be deduced from
one another by composition with a trivial deformation.}\end{defin}

Now, we get the following result:
\begin{prop}\label{deforminfieq} Two equivalent deformations have
$d''_A$-cohomologous infinitesimal deformation forms.\end{prop}
\begin{proof} For two equivalent deformations $B_t,\bar{B}_t$
there are local bases related as follows
$$\bar{b}_h(x,t)=exp(ta)(b_h(x,t)).$$ By applying the operator
$(\partial/\partial t)_{t=0}$ to this relation and using the
computation (\ref{auxeq}) we get the following relation between
the corresponding infinitesimal deformation forms
$$\bar{\Xi}^C=\Xi^C-d''_A(pr_Ca).$$ \end{proof}

The cohomology class $[\Xi]\in H^1(B;A/B)$ (the first cohomology
space of $A/B$-valued $B$-forms) is the {\it infinitesimal
deformation up to equivalence} of the foliation $B$ of the Lie
algebroid $A$.
\section{Integration of foliated Lie algebroids}
In this section we briefly discuss the integrability of a foliated
Lie algebroid.
\begin{defin}\label{grupoidfol} {\rm Let $l,r:G\rightrightarrows
M$ be a Lie groupoid ($l,r$ are the target (left) and source
(right) projections, respectively), where the unit manifold $M$ is
endowed with a foliation $\mathcal{F}$. We call $G$ an
$\mathcal{F}$-{\it foliated Lie groupoid} if it is endowed with a
foliation $
\mathcal{G}$ that has the following properties:

(i) the leaves of $ \mathcal{G}$ are contained in the left fibers,

(ii) $ \mathcal{G}$ is invariant by left translations,

(iii) the right projection
$r:(G,\mathcal{G})\rightarrow(M,\mathcal{F})$ is a foliated,
leaf-wise submersion.}\end{defin}

It is easy to get the following integrability result.
\begin{prop}\label{integralgfol} Let $l,r:(G,\mathcal{G})\rightrightarrows
(M,\mathcal{F})$ be a foliated Lie groupoid. Then, the Lie
algebroid $A(G)$ has a natural structure of a foliated Lie
algebroid. Conversely, if $(A,B)$ is a foliated Lie algebroid and
$A$ is integrable by the left-fiber connected Lie groupoid $G$,
then $G$ has a natural foliation $\mathcal{G}$ that makes $G$ into
a foliated Lie groupoid such that the restriction to the leaves of
$\mathcal{G}$ of the construction of $A(G)$ produces a Lie
subalgebroid that is isomorphic to $B$.\end{prop}
\begin{proof} If we start with the foliated groupoid $(G,\mathcal{G})$ and
construct the Lie algebroid $A(G)$ by means of the left-invariant
vector fields, the invariant vector fields that are tangent to the
leaves of $\mathcal{G}$ produce a Lie subalgebroid $B\subseteq
A(G)$.

Furthermore, if we denote by $T^lG$ the tangent bundle to the
$l$-fibers we get the exact sequence
$$0\rightarrow T\mathcal{G}\rightarrow T^lG\rightarrow T^lG/T\mathcal{G}
\rightarrow0,$$ where the quotient is a $\mathcal{G}$-foliated
bundle. The foliated cross sections of $T^lG$ are the
infinitesimal automorphisms of $T\mathcal{G}$ and (like for any
foliation) span the bundle $T^lG$. If the previous sequence is
quotientized by left translations we get the exact sequence
$$0\rightarrow B\rightarrow A(G)\rightarrow A(G)/B\rightarrow0$$
and we see that the pair $(A(G),B)$ satisfies the foliated
generation condition. Then, the anchor of $A(G)$ is induced by the
differential $r_*$ of the right projection and its restriction to
$B$ is onto $T\mathcal{F}$ because of condition (iii) of
Definition
\ref{grupoidfol}. Therefore, $(A(G),B)$ is a foliated Lie algebroid
over $(M,\mathcal{F})$.

Conversely, if we start with a foliated Lie algebroid $(A,B)$ over
$(M,\mathcal{F})$ and $A=A(G)$ for the Lie groupoid $G$, we get
the foliation $\mathcal{G}$ by asking it to be tangent to the left
invariant vector fields that produce cross sections of $B$ (see
Lemma 2.1 of
\cite{MM}). Then, the conditions required in Definition
\ref{grupoidfol} obviously hold. In particular, condition (iii)
follows from the foliated-generation condition of $(A,B)$ together
with the fact that $im\,\sharp_B=T\mathcal{F}$.\end{proof}
\section{Foliated Courant algebroids}
We do not intend to develop a theory of foliated Courant algebroid,
but, we would like to clarify how such a notion should be defined.
We refer to \cite{LWX} for the general theory of Courant algebroids.

Definition \ref{defalgtrans} of a transversal-Lie algebroid suggests
the following natural, analogous definition.
\begin{defin}\label{defCtransvers} {\rm A {\it transversal-Courant
algebroid} over the foliated manifold $(M,\mathcal{F})$ is a
foliated vector bundle $E\rightarrow M$ endowed with a symmetric,
non degenerate, foliated, inner product
$g\in\Gamma_{fol}\odot^2E^*$, with a foliated morphism
$\sharp_E:E\rightarrow\nu\mathcal{F}$ called the anchor and a
skew-symmetric bracket
$[\,,\,]_E:\Gamma_{fol}E\times\Gamma_{fol}E\rightarrow\Gamma_{fol}E$,
such that the following conditions are satisfied:

1)
$\sharp_E[e_1,e_2]_E=[\sharp_Ee_1,\sharp_Ee_2]_{\nu\mathcal{F}}$,

2) $im(\sharp_g\circ^t\sharp_E)\subseteq ker\,\sharp_E$,

3)
$\sum_{Cycl}[[e_1,e_2]_E,e_3]_E=(1/3)\partial\sum_{Cycl}g([e_1,e_2]_E,e_3)$,
$\partial f=(1/2)\sharp_g(^t\sharp_Edf)$,

4)
$[e_1,fe_2]_E=f[e_1,e_2]_E+((\sharp_Ee_1)f)e_2-g(e_1,e_2)\partial
f $,

5) $(\sharp_Ee)(g(e_1,e_2))=g([e,e_1]_E+\partial g(e,e_1)
,e_2)+g(e_1,[e,e_2]_E+\partial g(e,e_2)).$\\ In these conditions,
$e,e_1,e_2,e_3\in\Gamma_{fol}E$, $f\in
C^\infty_{fol}(M,\mathcal{F})$ and $t$ denotes
transposition.}\end{defin}

Using local, foliated, $g$-canonical bases, it is easy to see that
$^t\sharp_E(\lambda)\in\Gamma_{fol}E^*$,
$\forall\lambda\in\Gamma_{fol}(ann\,T\mathcal{F})$. In particular,
$\forall f\in C^\infty_{fol}(M,\mathcal{F})$, $\partial
f\in\Gamma_{fol}E$.

Like for the Courant algebroids, which are obtained if $\mathcal{F}$
is the foliation of $M$ by points, we will say that the
transversal-Courant algebroid $E$ is {\it transitive} if $\sharp_E$
is surjective and it is {\it exact} if it is transitive and
$rank\,E=2rank\,\nu\mathcal{F}$. The reason for this name is that
$E$ is exact if and only if the sequence of vector bundles
\begin{equation}\label{siralgexact} 0\longrightarrow ann\,T\mathcal{F}
\stackrel{\sharp_{(E,g)}}{\longrightarrow}E
\stackrel{\sharp_{E}}{\longrightarrow}\nu\mathcal{F}\longrightarrow0,\end{equation}
where $\sharp_{(E,g)}=\sharp_g\circ^t\hspace{-1pt}\sharp_E$
(therefore, $\partial f=\sharp_{(E,g)}df$) is exact. (This follows
like for Courant algebroids, e.g., Section 3 of \cite{V5}.)
\begin{rem}\label{Cexact} {\rm For an exact Courant algebroid any
choice of a splitting of the exact sequence (\ref{siralgexact}) that
has an isotropic image yields an isomorphism $E\approx(TM\oplus
T^*M)$ where the latter is endowed with a twisted Courant bracket
\cite{BCG}. The brackets of the transitive Courant algebroids were
determined in \cite{{S},{V5}} and may be made explicit by using a
splitting as above. These results extend to transversal Courant
algebroids if we add the hypothesis that the exact sequence
(\ref{siralgexact}) has a foliated splitting, which is no more an
automatic fact.}\end{rem}
\begin{rem}\label{obsproiectiiC} {\rm If $E$ is a transversal-Courant algebroid
over $(M,\mathcal{F})$, the local projected bundles $\{E_{Q_U}\}$
are Courant algebroids over the local transversal manifolds $Q_U$ of
$\mathcal{F}$. In the case of a submersion $M\rightarrow Q$ with
connected and simply connected fibers, $E\rightarrow M$ projects to
a Courant algebroid over $Q$.}\end{rem}
\begin{example}\label{CnuF} {\rm For any foliated manifold
$(M,\mathcal{F})$ the vector bundle
$\nu\mathcal{F}\oplus\nu^*\mathcal{F}=\nu\mathcal{F}\oplus
ann\,T\mathcal{F}$ has the structure of a transversal-Courant
algebroid induced by the classical Courant structure of the local
transversal manifolds $Q_U$ of Remark
\ref{obsproiectiiC}.}\end{example}
\begin{example}\label{extrCfol} {\rm Let $A$ be a Courant
algebroid over the foliated manifold $(M,\mathcal{F})$ and let $B$
be a subbundle of $A$ that is a Courant algebroid with respect to
the induced Courant bracket, anchor and metric. In particular, the
metric induced in $B$ must be non degenerate and $A=B\oplus C$,
$C=B^{\perp_g}$. A cross section $c\in\Gamma C$ will be called
$B$-foliated if, $\forall b\in\Gamma B$, $[b,c]_A\in\Gamma B$. By
the axioms of Courant algebroids (see condition 4) of Definition
\ref{defCtransvers}) we see that, on one hand, if the previous
condition holds for $b$ it also holds for $fb$ $(f\in C^\infty(M)$
and, on the other hand, if $c$ is $B$-foliated and $f\in
C^\infty_{fol}(M,\mathcal{F})$, $fc$ is foliated (use $g(b,c)=0$).
(The definition of a $B$-foliated cross section is not correct for
an arbitrary $a\in\Gamma A$.) Now, it makes sense to assume that
the pair $(A,B)$ satisfies the following conditions (analogous to
those in Definition
\ref{folLiealg}): i) $\sharp_A(B)=T\mathcal{F}$, ii) $\Gamma C$ is
locally spanned by the $B$-foliated cross sections. Like for Lie
algebroids (Section 1), these conditions show that $C$ has a
natural structure of a foliated vector bundle and we may add
condition iii) $g|_C\in\Gamma_{fol}\odot^2 C^*$. From i), we get
$^t
\sharp_A(ann\,T\mathcal{F})\subseteq ann\,B=C^*$, whence, $\partial
f\in\Gamma C$, $\forall f\in C^\infty_{fol}(M,\mathcal{F})$.
Furthermore, since $A$ is a Courant algebroid we have (see condition
5), Definition \ref{defCtransvers})
\begin{equation}\label{axv}
(\sharp_Ac_1)(g(b,c_2))=0=g([c_1,b]_A+\partial g(c_1,b),c_2)
+g(b,[c_1,c_2]_A+\partial g(c_1,c_2)),\end{equation} where
$b\in\Gamma B,c_1,c_2\in\Gamma C$ and $c_1,c_2$ are $B$-foliated.
In view of condition iii), (\ref{axv}) becomes
$g(b,[c_1,c_2]_A)=0$, therefore, $[c_1,c_2]_A\in\Gamma C$.
Moreover, using axiom 3) for the Courant algebroid $A$ and for the
triple $(c_1,c_2,b)$ instead of $(e_1,e_2,e_3)$ we see that the
bracket $[c_1,c_2]_A$ is $B$-foliated. Thus, conditions i), ii),
iii) imply that $(C,\sharp_A|_C,[\,,\,]_A|_{\Gamma_{fol}C})$ is a
transversal-Courant algebroid over
$(M,\mathcal{F})$.}\end{example}

We have formulated the above example because, at the first sight,
it could indicate the way to the notion of a foliated Courant
algebroid. However, it seems that there are no corresponding,
interesting, concrete examples and we shall propose another
procedure below.

Assume that $B$ is a $g$-isotropic subbundle of the Courant
algebroid $A$ of basis $(M,\mathcal{F})$ such that $\Gamma B$ is
closed by $A$-brackets and consider the $g$-coisotropic subbundle
$C=B^{\perp_g}\supseteq B$. By replacing $c_1$ by $b'\in\Gamma B$ in
(\ref{axv}), we see that $\Gamma B$ is closed by $A$-brackets if and
only if
\begin{equation}\label{echivalintegrab} [b,c]_A\in\Gamma C,\hspace{5mm}
\forall b\in\Gamma B,\forall c\in\Gamma C.\end{equation} Like in
Example \ref{extrCfol}, a cross section $c\in\Gamma C$ will be
called $B$-foliated if, $\forall b\in\Gamma B$, $[b,c]_A\in\Gamma B$
and this definition is correct. We will denote by $\Gamma_BC$ the
space of $B$-foliated cross sections of $C$.
\begin{defin}\label{Cfol} {\rm Let $A$ be a Courant algebroid over
the foliated manifold $(M,\mathcal{F})$. A subbundle $B\subseteq
A$ will be called a {\it foliation} of $A$ if the following
conditions are satisfied:

i) $B$ is a $g$-isotropic and $\Gamma B$ is closed by $A$-brackets,

ii) $\sharp_A(B)=T\mathcal{F}$,

iii) $\Gamma C$, $C=B^{\perp_g}$, is locally spanned by the
$B$-foliated cross sections,

iv) for any pair of $B$-foliated cross sections $c_1,c_2\in\Gamma C$
one has $g(c_1,c_2)\in C^\infty_{fol}(M,\mathcal{F})$.\\
\noindent A {\it foliated Courant algebroid} is a pair $(A,B)$ where
$A$ is a Courant algebroid and $B$ is a foliation of
$A$.}\end{defin}
\begin{prop}\label{CsitrC} If $(A,B)$ is a foliated Courant
algebroid over $(M,\mathcal{F})$, the vector bundle $E=C/B$,
$C=B^{\perp_g}$, inherits a natural structure of an $
\mathcal{F}$-transversal-Courant algebroid.\end{prop}
\begin{proof} Since $B$ is isotropic, there exist decompositions
\begin{equation}\label{descABB'} A=B\oplus S\oplus B',\;B\oplus
S=C,\;(B\oplus B') \perp_g S,\end{equation} where $B'$ is isotropic
and $g$ is non degenerate on the components $B\oplus B'$ and $S$. We
fix such a decomposition and define a structure of
transversal-Courant algebroid on $S$.

Take
$$g_S=g|_S,\;\sharp_S=\sharp_A\:(mod.\,T\mathcal{F}),\;[\,,\,]_S=
pr_S[\,,\,]_A,$$ where the projection on $S$ is defined by
(\ref{descABB'}). Like in the proof of Proposition \ref{ApeB},
properties ii) and iii) of Definition \ref{Cfol} yield a foliated
structure of $S$ such that $s\in\Gamma S$ is foliated with respect
to this structure if and only if $s$ is $B$-foliated. We check that,
$\forall s\in\Gamma_BC$, $\sharp_Ss$ is $ \mathcal{F}$-foliated.
Indeed, if $Y\in\Gamma T\mathcal{F}$, property ii) allows us to
write $Y=\sharp_Ab$, $b\in\Gamma B$ and we have
$$[Y,\sharp_As]_A=[\sharp_Ab,\sharp_As]_A=\sharp_A[b,s]_A
\in\Gamma T\mathcal{F},$$ therefore $\sharp_As$ is a foliated vector
field on $(M,\mathcal{F})$. Thus, $\sharp_S$ is a foliated morphism.
Furthermore, property iv) implies that $g_S$ is a foliated metric.

The decomposition (\ref{descABB'}) also gives $A^*=B^*\oplus
S^*\oplus B^{'*}$ and we have $$\flat_g(B)=B^{'*},\,\flat_g(S)=S^*,
\,\flat_g(B')=B^{*}.$$ If $\gamma\in ann\,T\mathcal{F}$,
$^t\sharp_A\gamma$ vanishes on $B$, because of property ii), hence,
$^t\sharp_A\gamma\in S^*\oplus B^{'*}$ and
$\sharp_g^t\sharp_A\gamma\in C$. In particular, $\forall f\in
C^\infty_{fol}(M,\mathcal{F}$, $\partial f\in\Gamma C$. This fact
allows us to use again (\ref{axv}), while taking
$s_1,s_2\in\Gamma_BC$ instead of $c_1,c_2$, and we see that
$[s_1,s_2]_A\in\Gamma C$. Now, if we write down the equality iii),
Definition \ref{defCtransvers} for the Courant algebroid $A$ and for
$e_1,e_2,e_3$ replaced by $s_1,s_2\in\Gamma_BC,b\in\Gamma B$, the
right hand side vanishes and we remain with $[[s_1,s_2]_A,b]_A=0$,
which shows that $[\,,\,]_S$ takes foliated cross sections to
foliated cross sections.

Thus, $S$ is endowed with all the structures required by Definition
\ref{defCtransvers} and it remains to check the conditions 1)-5). It
is easy to check that $\partial_S=\partial_A$ on
$C^\infty_{fol}(M,\mathcal{F})$. Accordingly, 1)-5) are implied by
the similar properties of the Courant algebroid $A$ (in particular,
for 2) we may use the fact that this condition is equivalent with
$g(\partial_Sf_1,\partial_Sf_2)=0$, $\forall f_1,f_2\in
C^\infty_{fol}(M,\mathcal{F})$ \cite{LWX}). The conclusion of the
proposition follows by transferring the structure obtained on $S$ to
$E$ via the natural isomorphism $E\approx S$.
\end{proof}
\begin{corol}\label{obsptC} If $(A,B)$ is a foliated Courant
algebroid over $(M,\mathcal{F})$ and if $\mathcal{F}$ consists of
the fibers of a submersion $M\rightarrow Q$ with connected and
simply connected fibers, the vector bundle $E=B^{\perp_g}/B$
projects to a Courant algebroid on $Q$.\end{corol}
\begin{rem}\label{obsprop} {\rm If $B$ is a subbundle of the
Courant algebroid $A$ that satisfies conditions i), ii) of
Definition \ref{Cfol} and if we ask the quotient bundle $A/B$ to
be a transversal-Courant algebroid such that its metric, anchor
and bracket be induced by those of $A$, then, obviously, $B$ is a
foliation of $A$.}\end{rem}

The previous proposition and corollary show the interest of the
notion of a foliated Courant algebroid. The following examples show
that Definition \ref{Cfol} is reasonable. In these examples
$T^{big}M=TM\oplus T^*M$ is the Courant algebroid defined in
\cite{Cou}; the anchor is the projection on $TM$ and the metric and
bracket are defined by \begin{equation}\label{gC}
g((X_1,\alpha_1),(X_2,\alpha_2))=\frac{1}{2}(\alpha_1(X_2)+\alpha_2(X_1)),
\end{equation}
\begin{equation}\label{crosetCou}
[(X_1,\alpha_1),(X_2,\alpha_2)]=([X_1,X_2], L_{X_1}\alpha_2-
L_{X_2}\alpha_1+\frac{1}{2}d(\alpha_1(X_2)-\alpha_2(X_1))).
\end{equation}
\begin{example}\label{exCfol0} {\rm A regular Dirac structure
$D\subseteq T^{big}M$ is a foliation of the Courant algebroid
$T^{big}M$ over $(M,\mathcal{E})$, where $\mathcal{E}$ is the
characteristic foliation of $D$. In this case the quotient bundle
is $D^{\perp_g}/D=0$.}\end{example}
\begin{example}\label{exCfol1} {\rm If $ \mathcal{F}$ is a foliation
of $M$, $T\mathcal{F}$ is a foliation of the Courant algebroid
$T^{big}M$. Indeed, using (\ref{gC}) we get
$C=T^{\perp_g}\mathcal{F}=TM\oplus ann\,T\mathcal{F}$. A pair
$(X,\gamma)\in\Gamma C$ is $T\mathcal{F}$-foliated if and only if
$X$ is an $ \mathcal{F}$-foliated vector field and $\gamma$ is an $
\mathcal{F}$-foliated $1$-form. The conditions of Definition
\ref{Cfol} are obviously satisfied and the corresponding
transversal-Courant algebroid is
$E=C/T\mathcal{F}=\nu\mathcal{F}\oplus ann\,T\mathcal{F}$, already
mentioned in Example \ref{CnuF}. Notice that $T\mathcal{F}$ remains
a foliation of $T^{big}M$ if the Courant bracket is twisted by means
of a closed $3$-form $\Phi$ such that $i(Y)\Phi=0$, $\forall Y\in
T\mathcal{F}$ \cite{SW}.}\end{example}
\begin{example}\label{exCfol2} {\rm Let $ \mathcal{F}$ be a foliation
of $M$ and $\theta\in\Omega^2(M)$ be a closed $2$-form. From
\cite{V1}, it is known that $E_\theta=\{(Y,i(Y)\theta)\,/\,Y\in
T\mathcal{F}\}$ is an isotropic subbundle of $T^{big}M$, which is
closed by the Courant bracket (\ref{crosetCou}), and its
orthogonal bracket is $E'_\theta=\{(X,i(X)\theta+\gamma)\,/\,X\in
TM,\gamma\in ann\,T\mathcal{F}\}$. A simple computation gives
$$[(Y,i(Y)\theta),(X,i(X)\theta+\gamma)]=([Y,X],i([Y,X])\theta+L_Y\gamma).$$
Accordingly, $(X,i(X)\theta+\gamma)\in\Gamma E'_\theta$ is
$E_\theta$-foliated if and only if
$X\in\Gamma_{fol}TM,\gamma\in\Omega^1_{fol}(M)$ and, like in Example
\ref{exCfol1}, the conditions of Definition \ref{Cfol} are satisfied
and $E_\theta$ is a foliation of $T^{big}M$. Again, the
corresponding quotient bundle is $E'_\theta/E_\theta\approx
\nu\mathcal{F}\oplus ann\,T\mathcal{F}$.}\end{example}

Below, we prove a result that has the flavor of a reduction and was
inspired by \cite{Zambon}. Let $B$ be a foliation of the Courant
algebroid $A$ over the foliated manifold $(M,\mathcal{F})$. Let $N$
be a submanifold of the base manifold $M$, and consider the
restricted vector bundles $A_N,B_N,C_N=B_N^{\perp_{g_A}}$.
\begin{prop}\label{reducereC} With the notation above, assume that:

(i) $N$ is transversal to and has a clean intersection with the
leaves of the foliation $\mathcal{F}$,

(ii) (the {\it reduction hypothesis}) $\sharp_A(C_N)\subseteq TN$.\\
Then, the vector bundle $E_N=C_N/B_N$ is a transversal-Courant
algebroid over $(N,\mathcal{F}'=\mathcal{F}\cap N)$. \end{prop}
\begin{proof} Since $B$ is a foliation of $A$, there exist local bases
$(b_h\in\Gamma B,a_u\in\Gamma_{B}C)$ of $\Gamma C$ and the induced
local bases of cross sections of $E=C/B$ have local transition
functions of the form $[\tilde a_u]_{{\rm
mod.}\,B}=\alpha_u^v[a_u]_{{\rm mod.}\,B}$ such that $\alpha_u^v$
are constant on the leaves of $\mathcal{F}$. Then, $\alpha_u^v$ are
constant on the leaves of $\mathcal{F}'$ too and the local bases
$[a_u|_N]_{{\rm mod.}\,B_N}$ of $\Gamma E_N$ define an
$\mathcal{F}'$-foliated structure on $E_N$. Using these bases, we
see that a cross section of $E_N$ is foliated with respect to the
$\mathcal{F}'$-foliated structure of $E_N$ if and only if, locally,
it is of the form $[c|_N]_{{\rm mod.}\,B_N}$ where $c$ is
$B$-foliated. Moreover, the metric induced on $E_N$ by the metric
$g$ of $A$ obviously is a foliated metric.

In view of the reduction hypothesis, the mapping
$$[c]_{{\rm mod.}\,B}\mapsto[\sharp_Ac]_{{\rm mod.}\,T\mathcal{F}}
\hspace{2mm}(c\in\Gamma C)$$
produces a vector bundle morphism
$\sharp_{E_N}:E_N\rightarrow\nu\mathcal{F}'$, which will be the
anchor of the required Courant structure on $E_N$. Like in the
proof of Proposition \ref{CsitrC}, we can see that this anchor is
a foliated morphism.

Now, we will show the existence of a bracket $[\,,\,]_{C_N}:\Gamma
C_N\times\Gamma C_N\rightarrow\Gamma E_N$ induced by the Courant
bracket of $A$. The definition of this bracket is
\begin{equation}\label{crosetCN} [c_1,c_2]_{C_N}(x)=[\tilde c_1,\tilde
c_2]_A(x)\;(mod.\,B_x),\hspace{3mm}x\in N,\end{equation} where
$\tilde c_1,\tilde c_2$ are arbitrary extensions of $c_1,c_2$ to
$\Gamma C$. In order to prove that the result does not depend on
the choice of the extensions it suffices to show that if $\tilde
c_2|_N=0$ then $[\tilde c_1,\tilde c_2]_A(x)\in B_x$. To see that,
take a local basis $(b_h\in\Gamma B,a_u\in\Gamma_BC)$ of $\Gamma
C$ and put $\tilde c_2= f^hb_h+k^ua_u$ where $f^h,k^u$ vanish on
$N$. Using the axioms of a Courant algebroid (see 4) of Definition
\ref{defCtransvers}), under the conditions above, we remain with
$$[\tilde c_1,\tilde c_2]_A(x)=-\sum_ug_x(\tilde c_1(x),a_u(x))\partial
k^u(x).$$ Furthermore, for any $c\in\Gamma C$ we have
$$g(c,\partial k^u)=\frac{1}{2}(\sharp_Ac)k^u=0$$ (because
$\sharp_Ac\in TN$ and $k^u|_N=0$), therefore, $\partial k^u(x)\in
B_x$, which ends the proof of the correctness of the definition of
the bracket (\ref{crosetCN}).

The bracket (\ref{crosetCN}) produces a bracket of foliated cross
sections of $E_N$ as follows. If $\gamma\in\Gamma_{fol}E_N$, then,
in a neighborhood of $x\in N$, we may write
$\gamma=\sum_u\varphi^u(x^a)[a_u|_N]_{{\rm mod.}\,B_N}$, which
implies that $\gamma$ may be represented by a cross sections of
$C_N$ that has the $B$-foliated, local extension
$\tilde\gamma=\sum_u\varphi^u(x^a)a_u$. If we put
\begin{equation}\label{crosetfolCN} [\gamma_1,\gamma_2]_{E_N}(x)=
[\tilde{\gamma}_1|_N,\tilde{\gamma_2}|_N]_{C_N}(x)=
[\tilde{\gamma}_1,\tilde{\gamma_2}]_A(x)\;(mod.\,B_x),\end{equation}
we get a well defined bracket on $\Gamma_{fol}E_N$.

Thus, on $E_N$ we have all the components required by the
definition of a transversal-Courant algebroid over $(N,
\mathcal{F}')$ and, moreover, these components may be seen as
restrictions to $N$ of the components of the transversal-Courant
algebroid $E$ over $(M,\mathcal{F})$ given by Proposition
\ref{CsitrC}. Hence, condition 1)-5) of Definition
\ref{defCtransvers} are satisfied and we are done.
\end{proof}

We finish by the observation that a definition similar to
Definition
\ref{defCtransvers} may be given for a notion of {\it holomorphic
Courant algebroid}; we just have to replace the word ``foliated"
by ``holomorphic" everywhere. Furthermore, if we state Definition
\ref{Cfol} for complex vector bundles $A,B$ over a complex analytic
manifold $M$, and with $T\mathcal{F}$ replaced by the
anti-holomorphic tangent bundle of $M$, we get the notion of a
foliation of a complex Courant algebroid. The corresponding
quotient $E$ will be a holomorphic Courant algebroid. Examples
\ref{exCfol1}, \ref{exCfol2} can be adapted to the
holomorphic case

\hspace*{7.5cm}{\small \begin{tabular}{l} Department of
Mathematics\\ University of Haifa, Israel\\ E-mail:
vaisman@math.haifa.ac.il \end{tabular}}
\end{document}